\documentclass{amsart}
\usepackage{amsfonts}
\usepackage{hyperref}
\usepackage{xr}

\newtheorem{theorem}{Theorem}[section]
\theoremstyle{plain}

\newtheorem{condition}{Condition}[section]

\newtheorem{corollary}{Corollary}[section]

\newtheorem{example}{Example}[section]

\newtheorem{lemma}{Lemma}[section]

\newtheorem{problem}{Problem}[section]
\newtheorem{proposition}{Proposition}[section]
\newtheorem{remark}{Remark}[section]

\numberwithin{equation}{section}

\begin{document}
\title[\textbf{Numbers of periodic orbits hidden at fixed points}]{%
\textbf{The numbers of periodic orbits of holomorphic mappings hidden at
fixed points}}
\author{Guang Yuan\medskip\ Zhang }
\address{Department of Mathematical Sciences, Tsinghua University, Beijing
100084, P. R. China. \textit{Email:} \textit{gyzhang@math.tsinghua.edu.cn}}

\begin{abstract}
Let $\Delta ^{2}$ be a ball in the complex vector space $\mathbb{C}^{2}$
centered at the origin, let $f:\Delta ^{2}\rightarrow \mathbb{C}^{2}$ be a
holomorphic mapping$,$ with $f(0)=0$, and let $M$ be a positive integer. If
the origin $0$ is an isolated fixed point of the $M$ th iteration $f^{M}$ of
$f,$ then one can define the number $\mathcal{O}_{M}(f,0)$ of periodic
orbits of $f$ with period $M$ hidden at the fixed point $0$, which has the
meaning: any holomorphic mapping $g:\Delta ^{2}\rightarrow \mathbb{C}^{2}$
sufficiently close to $f$ in a neighborhood of the origin has exactly $%
\mathcal{O}_{M}(f,0)$ distinct periodic orbits with period $M$ near the
origin, provided that all fixed points of $g^{M}$ near the origin are all
simple.

It is known that $\mathcal{O}_{M}(f,0)\geq 1$ iff the linear part of $f$ at
the origin has a periodic point of period $M.$ This paper will continue to
study the number $\mathcal{O}_{M}(f,0)$. We are interested in the condition
for the linear part of $f$ at the origin such that $\mathcal{O}_{M}(f,0)\geq
2.$

For a $2\times 2$ matrix $A$ that is arbitrarily given, the goal of this
paper is to give a necessary and sufficient condition for $A$,$\ $such that $%
\mathcal{O}_{M}(f,0)\geq 2$ for all holomorphic mappings $f:\Delta
^{2}\rightarrow \mathbb{C}^{2}$ such that $f(0)=0,$ $Df(0)=A$ and that the
origin $0$ is an isolated fixed point of $f^{M}.$
\end{abstract}

\thanks{Project 10271063 and 10571009 supported by NSFC}
\maketitle

\section{Introduction and the main Theorem\label{S1}}

We denote by $\mathbb{C}^{n}$ the complex vector space of dimension $n$ and
by $\mathcal{O}(\mathbb{C}^{n},0,0)$ the space of all germs of holomorphic
mappings $f$ from a neighborhood of the origin $0$ in $\mathbb{C}^{n}$ into $%
\mathbb{C}^{n}$ such that%
\begin{equation*}
f(0)=0\ \mathrm{(}0\mathrm{\ denotes\ the\ origin).}
\end{equation*}

Let%
\begin{equation*}
f(z)=\lambda z+\mathrm{higher\ terms},
\end{equation*}%
be a germ in $\mathcal{O}(\mathbb{C},0,0)$. Then $0$ is a fixed point of $f$
and for each $m\in \mathbb{N}$ (the set of positive integers)$,$ the $m$-th
iteration $f^{m}$ of $f$ is well defined in a neighborhood of $0$. $f^{m}$
is defined as $f^{1}=f,f^{2}=f\circ f,\dots ,f^{m}=f\circ f^{m-1},$
inductively.

If $\lambda =f^{\prime }(0)$ is a primitive $m$-th root of unity, then it is
well known in the theory of one variable complex dynamics [\ref{Mil}] that
either $f^{m}(z)\equiv z,$ or there exist an $\alpha \in \mathbb{N}$ and a
constant $a\neq 0$ such that
\begin{equation*}
f^{m}(z)=z+az^{\alpha m+1}+\mathrm{higher\ terms},
\end{equation*}%
and in the later case, $0$ is an isolated fixed point of $f^{m}$ and one can
split this fixed point into one fixed point and $\alpha $ periodic orbits by
small perturbations. More precisely, there exists a sequence of holomorphic
functions $f_{k}$ converging to $f$ uniformly in a neighborhood of the
origin so that $f_{k}(0)=0$ and each $f_{k}$ has $\alpha $ distinct periodic
orbits of period $m$ converging to the origin, and $\alpha $ is the largest
integer with this property. Note that $f$ itself has no periodic orbit of
period $m$ in a neighborhood of $0.$ The number $\alpha $ then can be
interpreted to be the number of periodic orbits of period $m$ of $f$
"hidden" at $0.$ Some authors call this phenomenon that $f$ has $\alpha $
\emph{virtual} periodic orbits of period $m$ at $0$ (see [\ref{CMY}])$.$
This number $\alpha $ will be denoted by $\mathcal{O}_{m}(f,0).$ Here, a
\emph{periodic orbit} of $f_{k}$ of period $m$ means a set $%
E=\{p_{1},p_{2},\dots ,p_{m}\}\subset B$ with cardinality $m$ such that $%
f_{k}(p_{1})=p_{2},f_{k}(p_{2})=p_{3},\dots ,f_{k}(p_{m-1})=p_{m}$ and $%
f_{k}(p_{m})=p_{1}.$ A point is called a \emph{periodic point of period} $m$
if and only if it is contained in a periodic orbit of period $m.$

In higher dimensional cases, there are similar phenomena, but things are far
more complicated. Now, let $f\in \mathcal{O}(\mathbb{C}^{n},0,0)$ and let $%
M\in \mathbb{N}$. Then there is a small ball $B$ centered at $0$ so that the
$M$-th iteration $f^{M}$ of $f$ is well defined in $\overline{B}$. If, in
addition, $0$ is an isolated fixed point of $f^{M}$ then we may make $B$
even smaller so that $0$ is the unique fixed point of $f^{M}$ in $\overline{B%
}$ and then define $\mathcal{O}_{M}(f,0)$ to be%
\begin{equation*}
\max \left\{ m;%
\begin{array}{l}
\text{\textrm{there\ exists\ a\ sequence }}f_{k}\in \mathcal{O}(\mathbb{C}%
^{n},0,0)\ \mathrm{uniformly\ converging\ to\ }f\ \text{\textrm{in}}\  \\
B\text{\textrm{\ }}\mathrm{such\ that\ each\ }f_{k}\ \mathrm{has\ }m\mathrm{%
\ }\text{\textrm{distinct\ }}\mathrm{periodic\ orbits\ of\ period\ }M\
\mathrm{in\ }B.%
\end{array}%
\right\}
\end{equation*}%
It is clear that in this definition, all periodic orbits of period $M$ of $%
f_{k}$ located in $B$ converge to $0$ uniformly as $k\rightarrow \infty .$
Therefore, this definition is independent of $B.$

The number $\mathcal{O}_{M}(f,0)$ is a well defined integer, and the
definition for $\mathcal{O}_{M}(f,0)$ agrees with that in the case $n=1.$ In
next section, we shall give another equivalent definition of the number $%
\mathcal{O}_{M}(f,0)$ and give some examples for understanding this number.

The following theorem is proved by the author in [\ref{Zh2}].

\begin{theorem}
\label{Th0.1}Let $f\in \mathcal{O}(\mathbb{C}^{n},0,0)$ and assume that the
origin is an isolated fixed point of $f^{M}.$ Then, $\mathcal{O}%
_{M}(f,0)\neq 0\ $if and only if the linear part of $f$ at $0$ has a
periodic point of period $M.$
\end{theorem}

The term \textquotedblleft linear part\textquotedblright\ indicates the
linear mapping $l:\mathbb{C}^{n}\rightarrow \mathbb{C}^{n},$
\begin{equation*}
l(x_{1},\dots ,x_{n})=(\sum_{j=1}^{n}a_{1j}x_{j},\dots
,\sum_{j=1}^{n}a_{nj}x_{j}),
\end{equation*}%
where
\begin{equation*}
(a_{ij})=Df(0)=\left. \left( \frac{\partial f_{i}}{\partial x_{j}}\right)
\right\vert _{0}
\end{equation*}%
is the Jacobian matrix of $f=(f_{1},\dots ,f_{n})$ at the origin$.$ If $M>1,$
then by Lemma \ref{j7}, the linear part of $f$ at $0$ has a periodic point
of period $M\ $if and only if the following condition holds.

\begin{condition}
\label{1.0} $Df(0)$ has eigenvalues $\lambda _{1},\dots ,\lambda _{s},s\leq
n,$ that are primitive $m_{1}$-th,..., $m_{s}$-th roots of unity,
respectively, such that $M$ is the least common multiple of $m_{1},\dots
,m_{s}.$
\end{condition}

Thus, if $M>1,$ by the above theorem, $\mathcal{O}_{M}(f,0)\geq 1$ if and
only if Condition \ref{1.0} holds. This gives rise to the following problem.

\begin{problem}
\label{P1}Assume that Condition \ref{1.0} holds. Under which additional
condition for $Df(0),\ $one must have $\mathcal{O}_{M}(f,0)\geq 2?$
\end{problem}

In this paper, we study this problem for the case $n=2$, and our goal is to
prove the following theorem.

\begin{theorem}
\label{Th1}Let $M>1$ be a positive integer and let $A$ be a $2$ by $2$
matrix. Then the following two conditions \textrm{(A) and (B) }are
equivalent:

\textrm{(A) }For any holomorphic mapping germ $f\in \mathcal{O}(\mathbb{C}%
^{2},0,0)$ such that $Df(0)=A$ and that $0$ is an isolated fixed point of
both $f$ and $f^{M},$
\begin{equation*}
\mathcal{O}_{M}(f,0)\geq 2.
\end{equation*}

\textrm{(B) }The two eigenvalues $\lambda _{1}$ and $\lambda _{2}$ of $A$
are primitive $m_{1}$ th and $m_{2}$ th roots of unity, respectively, and
one of the following conditions holds.

\textrm{(b1) }$A$ is diagonalizable, $m_{1}=m_{2}=M$ and $\lambda
_{1}=\lambda _{2}.$

\textrm{(b2) }$m_{1}=m_{2}=M$ and there exist positive integers $\alpha $
and $\beta $ such that $1<\alpha <M$, $1<\beta <M$ and%
\begin{equation*}
\lambda _{1}^{\alpha }=\lambda _{2},\lambda _{2}^{\beta }=\lambda
_{1},\alpha \beta >M+1.
\end{equation*}

\textrm{(b3) }$m_{1}|m_{2},$ $m_{2}=M$, and $\lambda _{2}^{m_{2}/m_{1}}\neq
\lambda _{1}.$

\textrm{(b4) }$M=[m_{1},m_{2}],(m_{1},m_{2})>1$ and $\max \{m_{1},m_{2}\}<M.$
\end{theorem}

Here, $[m_{1},m_{2}]$ denotes the least common multiple and $(m_{1},m_{2})$
denotes the greatest common divisor, of $m_{1}$ and $m_{2},$ and $%
m_{1}|m_{2} $ means that $m_{1}$ divides $m_{2}.$

\begin{remark}
\label{r1}Consider condition (b2). If in the theorem $m_{1}=m_{2}=M$ but $%
\lambda _{1}\neq \lambda _{2},$ say, $\lambda _{1}$ and $\lambda _{2}$ are
distinct primitive $M$ th roots of unity, then it is easy to show that there
uniquely exist positive integers $\alpha $ and $\beta $ such that $1<\alpha
<M$, $1<\beta <M$ and%
\begin{equation}
\lambda _{1}^{\alpha }=\lambda _{2},\lambda _{2}^{\beta }=\lambda _{1}.
\label{j6}
\end{equation}%
The proof is left to the reader.

(\ref{j6}) implies that $\lambda _{1}^{\alpha \beta }=\lambda _{1},$ and
then by the property of primitive $M$ th roots of unity, one can see that $%
\alpha \beta =kM+1$ for some positive integer $k$. Condition (b2) just
permits $k>1.$
\end{remark}

This paper is arranged as follows. In Section \ref{S2}, we shall give
another equivalent definition of the number $\mathcal{O}_{M}(f,0)$ and give
two examples for understanding this number. Sections \ref{S3}--\ref{PM} are
aimed to prove the main theorem. In Sections \ref{S3} and \ref{norm}, we
shall introduce some known results and prove a few consequences. Then, in
Section \ref{cronin-1}, we shall apply Cronin's theorem to compute zero
orders of some germs in $\mathcal{O}(\mathbb{C}^{2},0,0)$. After these
preparations, the proof of the main theorem will be given in the last two
sections.

\section{Another Definition of $\mathcal{O}_{M}(f,0)$ via Dold's Indices and
Some Examples\label{S2}}

Let $f\in \mathcal{O}(\mathbb{C}^{2},0,0)$. Then each component of $f$ can
be expressed as a power series at the origin. All power series in this paper
will be assumed to be convergent in a neighborhood of the origin.

If $p$ is an isolated zero of $f,$ say, there exists a ball $B$ centered at $%
p$ such that $f$ is well defined on $\overline{B}$ and that $p$ is the
unique solution of the equation $f(x)=0$ ($0$ is the origin) in $\overline{B}%
.$ Then we can define the \emph{zero order} (or \emph{multiplicity}) of $f$
at $p$ by%
\begin{equation*}
\pi _{f}(p)=\#\{x\in B;f(x)=q\},
\end{equation*}%
where $q$ is a regular value of $f$ such that $|q|$ is small enough and $\#$
denotes the cardinality. $\pi _{f}(p)$ is well defined (see \cite{LL} or [%
\ref{Zh}] for the details).

If the origin $0\ $is an isolated fixed point of $f,$ then $0$ is an
isolated zero of the germ $id-f\in \mathcal{O}(\mathbb{C}^{2},0,0),$ which
puts each $x$ near the origin into $x-f(x),$ and then the \emph{fixed point
index} of $f$ at $0$ is well defined by%
\begin{equation*}
\mu _{f}(0)=\pi _{id-f}(0).
\end{equation*}%
If $0$ is a fixed point of $f$ such that $id-f$ is regular at $0,$ say, the
Jacobian matrix $Df(0)$ of $f$ at $0$ has no eigenvalue $1,$ then $0$ is
called a \emph{simple} fixed point of $f.$ A simple fixed point of a
holomorphic mapping has index $1$ by the inverse function theorem (see Lemma %
\ref{lem2-a}). By the definition, it is clear that
\begin{equation*}
\mu _{f}(0)=\pi _{id-f}(0)=\pi _{f-id}(0).
\end{equation*}

If the origin $0$ is a fixed point of $f,$ then for any $m\in \mathbb{N},$
the $m$-th iteration $f^{m}$ is well defined in a neighborhood $V_{m}$ of $%
0. $ If for some $M\in \mathbb{N},$ the origin $0$ is an isolated fixed
point of both $f$ and $f^{M}$, then for each factor $m$ of $M,$ $0$ is an
isolated fixed point of $f^{m}$ as well and the fixed point index $\mu
_{f^{m}}(0)$ of $f^{m}$ at $0$ is well defined. Therefore, we can define the
(local) index as in A. Dold's work [\ref{Do}]:

\begin{equation}
P_{M}(f,0)=\sum_{\tau \subset P(M)}(-1)^{\#\tau }\mu _{f^{M:\tau }}(0),
\label{do-1}
\end{equation}%
where $P(M)$ is the set of all primes dividing $M,$ the sum extends over all
subsets $\tau $ of $P(M),$ $\#\tau \ $is the cardinal number of $\tau $ and $%
M:\tau =M(\prod_{p\in \tau }p)^{-1}$. Note that the sum includes the term $%
\mu _{f^{M}}(0)$ which corresponds to the empty subset $\tau =\emptyset $.
If $M=12=2^{2}\cdot 3,$ for example, then $P(M)=\{2,3\},$ and
\begin{equation*}
P_{12}(f,0)=\mu _{f^{12}}(0)-\mu _{f^{4}}(0)-\mu _{f^{6}}(0)+\mu _{f^{2}}(0).
\end{equation*}%
The formula (\ref{do-1}) is known as the M\"{o}bius inversion formula \cite%
{Hy} (see [\ref{Zh2}] for more interpretations of Dold's index).

\begin{remark}
(1) By Corollary \ref{do1},
\begin{equation*}
\mathcal{O}_{M}(f,0)=P_{M}(f,0)/M.
\end{equation*}

(2) By the definition, $\mathcal{O}_{1}(f,0)=P_{1}(f,0)=\mu _{f}(0).$
\end{remark}

We denote by $O(1)$ any holomorphic function germ at the origin, which may
be different in different places, even in a single equation; by $o(|z|^{k}),$
any holomorphic function germ $\alpha $ defined at the origin $z=0$ such
that
\begin{equation}
\lim_{z\rightarrow 0}|\alpha (z)|/|z|^{k}=0,  \label{1-1}
\end{equation}%
which is equivalent to the statement that $\alpha $ can be expressed as a
power series in which the terms of degree $\leq k$ are all zero. Also, the
same notation $o(|z|^{k})$ may denote different function germs in different
places, even in a single equation. When $o(|z|^{k})$ denotes a germ of one
variable function, we just write it to be $o(z^{k}).$ Thus, $o(1)=o(|z|^{0})$
means any holomorphic function vanishing at the origin.

\begin{example}
\label{e1}Let $f\in \mathcal{O}(\mathbb{C}^{2},0,0)$ be given by
\begin{equation*}
(x,y)\rightarrow (\lambda _{1}x+o(x),\lambda _{2}y+o(y)),
\end{equation*}%
such that $\lambda _{1}$ is a primitive $m_{1}$-th root of unity and $%
\lambda _{2}$ is a primitive $m_{2}$-th root of unity, $m_{1}$ and $m_{2}$
are distinct primes, and that $0$ is an isolated fixed point of the $%
m_{1}m_{2}$ th iteration $f^{m_{1}m_{2}}$ of $f.$ Then there exist nonzero
constants $a,b$ and positive integers $\alpha $ and $\beta $ such that%
\begin{equation*}
f^{m_{1}}(x,y)^{T}=\left(
\begin{array}{l}
x+ax^{\alpha m_{1}+1}(1+o(1)) \\
\lambda _{2}^{m_{1}}y(1+o(1))%
\end{array}%
\right) ,
\end{equation*}%
\begin{equation*}
f^{m_{2}}(x,y)^{T}=\left(
\begin{array}{l}
\lambda _{1}^{m_{2}}x(1+o(1)) \\
y+by^{\beta m_{2}+1}(1+o(1))%
\end{array}%
\right) ,
\end{equation*}%
and%
\begin{equation*}
f^{m_{1}m_{2}}(x,y)^{T}=\left(
\begin{array}{l}
x+am_{2}x^{\alpha m_{1}+1}(1+o(1)) \\
y+bm_{1}x^{\beta m_{2}+1}(1+o(1))%
\end{array}%
\right) .
\end{equation*}%
Thus, considering that $\lambda _{1}^{m_{2}}\neq 1$ and $\lambda
_{2}^{m_{1}}\neq 1$, we have by Cronin's Theorem introduced in Section \ref%
{cronin-1} that%
\begin{equation*}
\mu _{f}(0)=1,
\end{equation*}%
\begin{equation*}
\mu _{f^{m_{1}}}(0)=\alpha m_{1}+1,\ \mu _{f^{m_{2}}}(0)=\beta m_{2}+1,
\end{equation*}%
\begin{equation*}
\mu _{f^{m_{1}m_{2}}}(0)=\left( \alpha m_{1}+1\right) (\beta m_{2}+1),
\end{equation*}%
and then by the formula (\ref{do-1}),%
\begin{equation*}
P_{m_{1}}(f,0)=\mu _{f^{m_{1}}}(0)-\mu _{f}(0)=\alpha m_{1},
\end{equation*}%
\begin{equation*}
P_{m_{2}}(f,0)=\mu _{f^{m_{2}}}(0)-\mu _{f}(0)=\beta m_{2},
\end{equation*}%
\begin{equation*}
P_{m_{1}m_{2}}(f,0)=\mu _{f^{m_{1}m_{2}}}(0)-\mu _{f^{m_{1}}}(0)-\mu
_{f^{m_{2}}}(0)+1=\alpha \beta m_{1}m_{2},
\end{equation*}%
and then, by Corollary \ref{do1},
\begin{equation*}
\mathcal{O}_{m_{1}}(f,0)=\alpha ,\mathcal{O}_{m_{2}}(f,0)=\beta ,\mathcal{O}%
_{m_{1}m_{2}}(f,0)=\alpha \beta ,
\end{equation*}%
say%
\begin{equation}
\mathcal{O}_{m_{1}m_{2}}(f,0)=\mathcal{O}_{m_{1}}(f,0)\mathcal{O}%
_{m_{2}}(f,0).  \label{j1}
\end{equation}
\end{example}

By this example, one may guess that there is a relation between the numbers $%
\mathcal{O}_{m_{1}}(f,0),$ $\mathcal{O}_{m_{2}}(f,0)$ and $\mathcal{O}%
_{m_{1}m_{2}}(f,0)\ $similar to the above equality (\ref{j1}). But see the
next example.

\begin{example}
\label{E2}Let $k>1$ be any given positive integer and let $f\in \mathcal{O}(%
\mathbb{C}^{2},0,0)$ be given by%
\begin{equation*}
\left( f(x,y)\right) ^{T}=\left(
\begin{array}{c}
-x+x^{2k+1}+xy^{3} \\
e^{\frac{2\pi i}{3}}y+x^{2}y+y^{3k+1}%
\end{array}%
\right) .
\end{equation*}%
We show that $\mathcal{O}_{2}(f,0)=\mathcal{O}_{3}(f,0)=k,$ but $\mathcal{O}%
_{6}(f,0)=1.$

After a careful computation, we have%
\begin{equation}
\left( f^{2}(x,y)\right) ^{T}=\left(
\begin{array}{l}
x-2x^{2k+1}(1+o(1))-2xy^{3}(1+o(1)) \\
e^{\frac{4\pi i}{3}}y(1+o(1))%
\end{array}%
\right) ,  \label{ggg1}
\end{equation}%
\begin{equation}
\left( f^{3}(x,y)\right) ^{T}=\left(
\begin{array}{l}
-x(1+o(1)) \\
y+3e^{\frac{4\pi i}{3}}x^{2}y(1+o(1))+3e^{\frac{4\pi i}{3}}y^{3k+1}(1+o(1))%
\end{array}%
\right) ,  \label{ggg2}
\end{equation}%
and%
\begin{equation*}
\left( f^{6}(x,y)\right) ^{T}=\left(
\begin{array}{l}
x+xh_{1}(x,y) \\
y+yh_{2}(x,y)%
\end{array}%
\right) ,
\end{equation*}%
with%
\begin{equation*}
\begin{array}{l}
h_{1}(x,y)=-6x^{2k}(1+o(1))-6y^{3}(1+o(1)), \\
h_{2}(x,y)=6e^{\frac{4\pi i}{3}}x^{2}(1+o(1))+6e^{\frac{4\pi i}{3}%
}y^{3k}(1+o(1)).%
\end{array}%
\end{equation*}

We first compute $\mu _{f^{2}}(0).$ It equals the zero order of the mapping
\begin{equation*}
F(x,y)=(x,y)-f^{2}(x,y),
\end{equation*}%
and then, by (\ref{ggg1}) and Corollary \ref{c2}, we have
\begin{equation*}
\mu _{f^{2}}(0)=\pi _{F}(0)=2k+1.
\end{equation*}%
Similarly, by (\ref{ggg2}) and Corollary \ref{c2}, we have
\begin{equation*}
\mu _{f^{3}}(0)=3k+1.
\end{equation*}

On the other hand, $0$ is a simple fixed point of $f,$ and then by Lemma \ref%
{lem2-a}, $\mu _{f}(0)=1.$ Therefore, by the formula (\ref{do-1}), we have%
\begin{eqnarray*}
P_{2}(f,0) &=&\mu _{f^{2}}(0)-1=2k, \\
P_{3}(f,0) &=&\mu _{f^{3}}(0)-1=3k,
\end{eqnarray*}%
and then, by Corollary \ref{do1}, we have
\begin{eqnarray*}
\mathcal{O}_{2}(f,0) &=&P_{2}(f,0)/2=k, \\
\mathcal{O}_{3}(f,0) &=&P_{3}(f,0)/3=k.
\end{eqnarray*}

Next, we show that $\mathcal{O}_{6}(f,0)=1.$ It is clear that $\mu
_{f^{6}}(0)$ equals the zero order of the mapping
\begin{equation*}
id-f^{6}:(x,y)\mapsto -(xh_{1}(x,y),yh_{2}(x,y)),
\end{equation*}%
and by Lemma \ref{prod1}, the zero order of $id-f^{6}$ at $0$ is the sum of
the zero orders of the four mappings putting $(x,y)\ $into $(x,y),$ $%
(x,h_{2}(x,y)),$ $(h_{1}(x,y),y)$ and $(h_{1}(x,y),h_{2}(x,y)),$ which are $%
1,$ $3k$, $2k$ and $6,$ respectively, by Cronin's Theorem and Corollary \ref%
{c2}. Thus $\mu _{f^{6}}(0)=5k+7,$ and then, by the formula (\ref{do-1}), we
have%
\begin{eqnarray*}
P_{6}(f,0) &=&\mu _{f^{6}}(0)-\mu _{f^{2}}(0)-\mu _{f^{3}}(0)+1 \\
&=&5k+7-2k-1-3k-1+1,
\end{eqnarray*}%
and then $P_{6}(f,0)=6,$ and $\mathcal{O}_{6}(f,0)=1.$
\end{example}

\begin{remark}
(1). Assume that $m_{1}$ and $m_{2}$ are two positive integers and let $f\in
\mathcal{O}(\mathbb{C}^{2},0,0)$ such that 0 is an isolated fixed point of $%
f^{[m_{1},m_{2}]}$, where $[m_{1},m_{2}]$ denotes the least common multiple
of $m_{1}$ and $m_{2}$. Then, by Theorem \ref{Th0.1}, $\mathcal{O}%
_{m_{1}}(f,0)\geq 1$ and $\mathcal{O}_{m_{2}}(f,0)\geq 1$ imply $\mathcal{O}%
_{[m_{1},\ m_{2}]}(f,0)\geq 1$.

(2). By the main theorem, when $m_{1}$ and $m_{2}$ satisfy certain condition
(for example, if $m_{1}=6$ and $m_{2}=10),$ $\mathcal{O}_{m_{1}}(f,0)\geq 1$
and $\mathcal{O}_{m_{2}}(f,0)\geq 1$ implies $\mathcal{O}_{[m_{1},\
m_{2}]}(f,0)\geq 2.$
\end{remark}

\section{Some basic results of fixed point indices and zero orders\label{S3}}

In this section we introduce some results for later use. Most of them are
known.

Let $U$ be an open and bounded subset of $\mathbb{C}^{2}$ and let $f:%
\overline{U}\rightarrow \mathbb{C}^{2}$ be a holomorphic mapping. If $f$ has
no fixed point on the boundary $\partial U$, then the fixed point set
\textrm{Fix}$(f)$ of $f$ is a compact analytic subset of $U,$ and then it is
finite (see [\ref{Ch}]); and therefore, we can define the global fixed point
index $L(f)$ of $f$ as:
\begin{equation*}
L(f)=\sum_{p\in \mathrm{Fix}(f)}\mu _{f}(p),
\end{equation*}%
which is just the number of all fixed points of $f$, counting indices. $L(f)$
is, in fact, the Lefschetz fixed point index of $f$ (see the appendix
section in [\ref{Zh2}] for the details).

For each $m\in \mathbb{N},$ the $m$-th iteration $f^{m}$ of $f$ is
understood to be defined on
\begin{equation*}
K_{m}(f)=\cap _{k=0}^{m-1}f^{-k}(\overline{U})=\{x\in \overline{U}%
;f^{k}(x)\in \overline{U}\ \mathrm{for\ all\ }k=1,\dots ,m-1\},
\end{equation*}%
which is the largest set where $f^{m}$ is well defined. Since $U$ is
bounded, $K_{m}(f)$ is a compact subset of $\overline{U}$. Here, $f^{0}=id.$

Now, let us introduce the global Dold's index. Let $M\in \mathbb{N}$ and
assume that $f^{M}$ has no fixed point on the boundary $\partial U.$ Then,
for each factor $m$ of $M,$ $f^{m}$ again has no fixed point on $\partial U,$
and then the fixed point set $\mathrm{Fix}(f^{m})$ of $f^{m}$ is a compact
subset of $U$. Thus, there exists an open subset $V_{m}$ of $U$ such that $%
\mathrm{Fix}(f^{m})\subset V_{m}\subset \overline{V_{m}}\subset U$ and $%
f^{m} $ is well defined on $\overline{V_{m}},$ and thus $L(f^{m}|_{\overline{%
V_{m}}})$ is well defined and we write $L(f^{m})=L(f^{m}|_{\overline{V_{m}}%
}),$ where $f^{m}|_{\overline{V_{m}}}$ is the restriction of $f^{m}$ to $%
\overline{V_{m}}.$ In this way, we can define the global Dold's index (see [%
\ref{Do}]) as (\ref{do-1}):%
\begin{equation}
P_{M}(f)=\sum_{\tau \subset P(M)}(-1)^{\#\tau }L(f^{M:\tau }).  \label{p0p}
\end{equation}

Let $m\in \mathbb{N}$. It is clear that, for any compact subset $K$ of $U$
with $\cup _{j=1}^{m}f^{j}(K)\subset U$, there is a neighborhood $V\subset U$
of $K,$ such that for any holomorphic mapping $g:\overline{U}\rightarrow
\mathbb{C}^{2}$ sufficiently close to $f,$ the iterations $g^{j},j=1,\dots
,m,$ are well defined on $\overline{V}$ and%
\begin{equation*}
\max_{X\in \overline{U}}|g(X)-f(X)|\rightarrow 0\Longrightarrow \max_{1\leq
j\leq m}\max_{X\in \overline{V}}|g^{j}(X)-f^{j}(X)|\rightarrow 0.
\end{equation*}%
We shall use these facts frequently and tacitly.

We denote by $\Delta ^{2}$ a ball in $\mathbb{C}^{2}$ centered at the origin.

\begin{lemma}[\protect\cite{LL}]
\label{lem2-a}Let $f\in \mathcal{O}(\mathbb{C}^{2},0,0)$ and assume that the
origin is an isolated fixed point$.$ Then%
\begin{equation*}
\mu _{f}(0)\geq 1,
\end{equation*}%
and the equality holds if and only if $1$ is not an eigenvalue of $Df(0)$.
\end{lemma}

\begin{lemma}[\protect\cite{LL}]
\label{lem2-1}(1) Let $f:\overline{\Delta ^{2}}\rightarrow \mathbb{C}^{2}$
be a holomorphic mapping such that $f$ has no fixed point on the boundary $%
\partial \Delta ^{2}.$ Then there exists a $\delta >0$ such that any
holomorphic mapping $g:\overline{\Delta ^{2}}\rightarrow \mathbb{C}^{2}$
with $\max_{x\in \overline{\Delta ^{2}}}|g(x)-f(x)|<\delta $ has finitely
many fixed points in $\Delta ^{2}$ and satisfies
\begin{equation*}
L(g)=\sum_{p\in \mathrm{Fix}(g)}\mu _{g}(p)=\sum_{p\in \mathrm{Fix}(f)}\mu
_{f}(p)=L(f).
\end{equation*}

(2) In particular, if $0$ is the unique fixed point of $f$ in $\overline{%
\Delta ^{2}},$ then for any holomorphic mapping $g:\overline{\Delta ^{2}}%
\rightarrow \mathbb{C}^{2}$ with $\max_{x\in \overline{\Delta ^{2}}%
}|g(x)-f(x)|<\delta ,\ $%
\begin{equation*}
\mu _{f}(0)=\sum_{p\in \mathrm{Fix}(g)}\mu _{g}(p),
\end{equation*}%
and if in addition all fixed points of $g$ are simple, then
\begin{equation*}
\mu _{f}(0)=\#\mathrm{Fix}(g)=\#\{y\in \Delta ^{2};g(y)=y\}.
\end{equation*}
\end{lemma}

This result is another version of Rouch\'{e}'s theorem which is stated as
follows.

\begin{lemma}[Rouch\'{e}'s theorem \protect\cite{LL}]
Let $f:\overline{\Delta ^{2}}\rightarrow $ $\mathbb{C}^{2}$\ be a
holomorphic mapping such that $f$ has no zero on $\partial \Delta ^{2}.$
Then there exists a $\delta >0$ such that any holomorphic mapping $g:%
\overline{\Delta ^{2}}\rightarrow $ $\mathbb{C}^{2}$ with $\max_{x\in
\partial \Delta ^{2}}|g(x)-f(x)|<\delta $ has the same number of zeros in $%
\Delta ^{2}$ as $f,$ counting zero orders, say,%
\begin{equation*}
\sum_{f(x)=0}\pi _{f}(x)=\sum_{g(x)=0}\pi _{g}(x).
\end{equation*}
\end{lemma}

\begin{corollary}
Let $U$ be a bounded open subset of $\mathbb{C}^{2},$ let $M\in \mathbb{N},$
let $f:\overline{U}\rightarrow \mathbb{C}^{2}$ be a holomorphic mapping and
assume that $f^{M}$ has no fixed point on $\partial U.$ If $f$ has a
periodic point $p\in U$ with period $M$, then any holomorphic mapping $g:%
\overline{U}\rightarrow \mathbb{C}^{2}$ that is sufficiently close to $f$
has a periodic point with period $M$ in $U.$
\end{corollary}

\begin{proof}
This follows from Lemmas \ref{lem2-a} and \ref{lem2-1} directly (see [\ref%
{Zh2}] for a simple proof).
\end{proof}

\begin{corollary}
\label{adc}Let $U$ be a bounded open subset of $\mathbb{C}^{2},$ let $M\in
\mathbb{N},$ let $f:\overline{U}\rightarrow \mathbb{C}^{2}$ be a holomorphic
mapping and assume that $f^{M}$ has no fixed point on $\partial U.$ If $f$
has $k$ distinct periodic points with period $M$, then any holomorphic
mapping $g:\overline{U}\rightarrow \mathbb{C}^{2}$ that is sufficiently
close to $f$ has at least $k$ distinct periodic points with period $M$ in $%
U. $
\end{corollary}

\begin{proof}
This follows from the previous corollary directly.
\end{proof}

\begin{corollary}
\label{ad-2-0}Let $\{f_{n}\}\subset \mathcal{O}(\mathbb{C}^{2},0,0)$ be a
sequence converging to $f\in \mathcal{O}(\mathbb{C}^{2},0,0),$ uniformly in
a neighborhood of the origin. If the origin is an isolated zero of $f$, and
there exists an integer $k$ such that $\pi _{f_{n}}(0)\geq k$ for all $n\in
\mathbb{N},$ then $\pi _{f}(0)\geq k.$
\end{corollary}

\begin{proof}
This follows from Rouch\'{e}'s theorem directly.
\end{proof}

\begin{corollary}[{[\protect\ref{Zh2}]}]
\label{cc2-1}Let $M$ be a positive integer, let $U$ be a bounded open subset
of $\mathbb{C}^{2}$, let $f:\overline{U}\rightarrow \mathbb{C}^{2}$ be a
holomorphic mapping and assume that $f^{M}$ has no fixed point on $\partial
U $. Then:

(i). There exists an open subset $V$ of $U,$ such that $f^{M}$ is well
defined on $\overline{V},$ has no fixed point outside $V,$ and has only
finitely many fixed points in $V$.

(ii). For any holomorphic mapping $g:\overline{U}\rightarrow \mathbb{C}^{2}$
sufficiently close to $f$, $g^{M}$ is well defined on $\overline{V},$ has no
fixed point outside $V$ and has only finitely many fixed points in $V;$ and
furthermore,%
\begin{equation*}
L(g^{M})=L(f^{M})\mathrm{,}\text{ \ }P_{M}(g)=P_{M}(f).
\end{equation*}

(iii). In particular, if $p_{0}\in U$ is the unique fixed point of both $f$
and $f^{M}$ in $\overline{U},$ then for any holomorphic mapping $g:\overline{%
U}\rightarrow \mathbb{C}^{2}$ sufficiently close to $f$,
\begin{equation*}
L(g^{M})=L(f^{M})=\mu _{f^{M}}(p_{0})\mathrm{,}\text{ }%
P_{M}(g)=P_{M}(f)=P_{M}(f,p_{0}).
\end{equation*}
\end{corollary}

\begin{remark}
Under the assumption that $f^{M}$ has no fixed point on $\partial U,$ for
any factor $m$ of $M,$ the conclusions (i)--(iii) remain valid if $M$ is
replaced by $m,$ since $f^{m}$ has no fixed point on $\partial U$ as well.
\end{remark}

\begin{lemma}
\label{lem2-1-2}Let $M$ be a positive integer and let $f:\overline{\Delta
^{2}}\rightarrow \mathbb{C}^{2}$\ be a holomorphic mapping such that $f^{M}$
has no fixed point on $\partial \Delta ^{2}$ and each fixed point of $f^{M}$
is simple. Then, $\mathrm{Fix}(f^{M})$ is finite, and

\textrm{(i) }$L(f^{M})=\#\mathrm{Fix}(f^{M})=\sum_{m|M}P_{m}(f);$

\textrm{(ii) }$P_{M}(f)$ is the cardinal number of the set of periodic
points of $f$ of period $M$;

\textrm{(iii) }$P_{M}(f)/M$ is the number of distinct periodic orbits of $f$
of period $M.$
\end{lemma}

\begin{proof}
(i) and (ii) are proved in [\ref{FL}] (see [\ref{Zh2}] for a very simple
proof), and (iii) follows from (ii).
\end{proof}

A fixed point $p$ of $f$ is called \textit{hyperbolic} if $Df(p)$ has no
eigenvalue of absolute $1.$ If $p$ is a hyperbolic fixed point of $f,$ then
it is a hyperbolic fixed point of all iterations $f^{j},j\in \mathbb{N}.$ A
hyperbolic fixed point is a simple fixed point, and so it has index $1$ by
Lemma \ref{lem2-a}$.$

\begin{lemma}
\label{lem2-1-1}Let $M$ be a positive integer, let $V$ be an open subset of $%
\Delta ^{2}$ and let $f:$ $\overline{\Delta ^{2}}\rightarrow \mathbb{C}^{2}$
be a holomorphic mapping such that $f^{M}$ is well defined in $\overline{V}$
and has no fixed point on $\partial V$. Then for any positive number $%
\varepsilon >0,$ there exists a holomorphic mapping $f_{\varepsilon }:%
\overline{\Delta ^{2}}\rightarrow \mathbb{C}^{2}$ such that
\begin{equation*}
\sup_{X\in \overline{\Delta ^{2}}}|f_{\varepsilon {}}(X)-f(X)|<\varepsilon
\end{equation*}%
and all the fixed points of $f_{\varepsilon }^{M}$ located in $\overline{V}$
are hyperbolic.
\end{lemma}

A proof of this result follows from the argument in [\ref{B}]. Another proof
can be found in [\ref{Zh1}].

\begin{corollary}
\label{do1}Let $f:\overline{\Delta ^{2}}\rightarrow \mathbb{C}^{2}$\ be a
holomorphic mapping so that $0$ is the unique fixed point of both $f$ and $%
f^{M}$ in $\overline{\Delta ^{2}}.$ Then%
\begin{equation}
\mathcal{O}_{M}(f,0)=P_{M}(f,0)/M,  \label{j3}
\end{equation}%
and there exists a $\delta >0,$ such that any holomorphic mapping $g:%
\overline{\Delta ^{2}}\rightarrow \mathbb{C}^{2}$ with $\max_{x\in \overline{%
\Delta ^{2}}}|g(x)-f(x)|<\delta $ has exactly $\mathcal{O}_{M}(f,0)$
distinct periodic orbits of period $M$ in $\Delta ^{2},$ provided that all
fixed points of $g^{M}$ are simple.
\end{corollary}

\begin{proof}
By Corollary \ref{cc2-1} (iii), there exists a $\delta >0$ such that, for
any holomorphic mapping $g:\overline{\Delta ^{2}}\rightarrow \mathbb{C}^{2}$%
,
\begin{equation}
\max_{x\in \overline{\Delta ^{2}}}|g(x)-f(x)|<\delta  \label{d1}
\end{equation}%
implies
\begin{equation}
P_{M}(g)=P_{M}(f,0).  \label{d2}
\end{equation}

Let $\varepsilon $ be any positive number with $\varepsilon <\delta .$ Then
by Lemma \ref{lem2-1-1}, there exists a holomorphic mapping $g_{1}:\overline{%
\Delta ^{2}}\rightarrow \mathbb{C}^{2}$ satisfying (\ref{d1}) for $%
\varepsilon $, say,
\begin{equation*}
\max_{x\in \overline{\Delta ^{2}}}|g_{1}(x)-f(x)|<\varepsilon ,
\end{equation*}%
such that all fixed points of $g_{1}^{M}$ located in $\Delta ^{2}$ are
simple, and then by (\ref{d2}) and Lemma \ref{lem2-1-2} (iii), $g=g_{1}$ has
exactly $P_{M}(f,0)/M$ distinct periodic orbits of period $M.$ Therefore, by
the definition of the number $\mathcal{O}_{M}(f,0)$ and the arbitrariness of
$\varepsilon ,$ we have
\begin{equation*}
\mathcal{O}_{M}(f,0)\geq P_{M}(f,0)/M.
\end{equation*}%
We show that the inequality does not occur.

Otherwise, by the definition of $\mathcal{O}_{M}(f,0),$ there exists a
holomorphic mapping $g_{2}:\overline{\Delta ^{2}}\rightarrow \mathbb{C}^{2}$
satisfying (\ref{d1}), say,
\begin{equation*}
\max_{x\in \overline{\Delta ^{2}}}|g_{2}(x)-f(x)|<\delta ,
\end{equation*}%
such that $g_{2}$ has at least $P_{M}(f,0)+M$ distinct periodic points of
period $M$ in $\Delta ^{2}.$

Furthermore, by Corollary \ref{adc} and Lemma \ref{lem2-1-1}, there exists a
holomorphic mapping $g_{3}$ satisfying the following three conditions.

(a) $g_{3}$ is so close to $g_{2}$ that $g_{3}$ satisfies (\ref{d1}), say,
\begin{equation*}
\max_{x\in \overline{\Delta ^{2}}}|g_{3}(x)-f(x)|<\delta .
\end{equation*}

(b) $g_{3}$ is so close to $g_{2}$ that it has at least $P_{M}(f,0)+M$
distinct periodic points of period $M$ in $\Delta ^{2}.$

(c) All fixed points of $g_{3}^{M}$ in $\Delta ^{2}$ are simple.

By (c) and Lemma \ref{lem2-1-2} (ii), $g_{3}$ has exactly $P_{M}(g_{3})$
distinct periodic points of period $M,$ and then by (b), $P_{M}(g_{3})\geq
P_{M}(f,0)+M.$ But $g=g_{3}$ satisfies (\ref{d1}) by (a). Then $g=g_{3}$
satisfies (\ref{d2}), say, $P_{M}(g_{3})=P_{M}(f,0)$. This is a
contradiction, and then (\ref{j3}) is proved.

If $g:\overline{\Delta ^{2}}\rightarrow \mathbb{C}^{2}$ is a holomorphic
mapping that satisfies (\ref{d1}), then it satisfies (\ref{d2}), and then,
by Lemma \ref{lem2-1-2} (iii), it has
\begin{equation*}
P_{M}(g)/M=P_{M}(f,0)/M=\mathcal{O}_{M}(f,0)
\end{equation*}%
distinct periodic orbits in $\Delta ^{2}$ of period $M,$ provided that all
fixed points of $g^{M}$ are simple. This completes the proof.
\end{proof}

The following result also follows from the above argument.

\begin{lemma}
\label{lem2-2}Let $k$ and $M$ be positive integers and let $f\in \mathcal{O}(%
\mathbb{C}^{2},0,0)$. If $0$ is an isolated fixed point of both $f$ and $%
f^{M},$ then $\mathcal{O}_{M}(f,0)\geq k$ if there exists a sequence of
holomorphic mappings $f_{j}\in \mathcal{O}(\mathbb{C}^{2},0,0),$ uniformly
converging to $f$ in a neighborhood of the origin, such that
\begin{equation*}
f_{j}(0)=0\;\mathrm{and\;}\mathcal{O}_{M}(f_{j},0)\geq k.
\end{equation*}
\end{lemma}

\begin{lemma}
\label{j7}Let $L:\mathbb{C}^{n}\rightarrow \mathbb{C}^{n}$ be a linear
mapping and let $M>1\ $be a positive integer. Then $L$ has a periodic point
of period $M$ if and only if $L$ has eigenvalues $\lambda _{1},\dots
,\lambda _{s},s\leq n,$ that are primitive $m_{1}$ th,..., $m_{s}$ th roots
of unity, respectively, such that $M=\left[ m_{1},\dots ,m_{s}\right] .$
\end{lemma}

This is a basic knowledge of elementary linear algebra. Recall that $%
[m_{1},\dots ,m_{s}]$ denotes the least common multiple of $m_{1},\dots
,m_{s}$. This Lemma is only used once in this paper (in Section 1). We shall
frequently use its special case with $n=2$ in another version:

\begin{lemma}
\label{lem2-5}Let $L:\mathbb{C}^{2}\rightarrow \mathbb{C}^{2}$ be a linear
mapping and let $M>1\ $be a positive integer. Then $L$ has a periodic point
of period $M$ if and only if one of the following conditions holds.

$\mathrm{(a)}\;$One eigenvalue of $L$\ is a primitive $M$-th root of unity.

$\mathrm{(b)}$ The two eigenvalues of $L$\ are primitive $m_{1}$ th and $%
m_{2}$ th roots of unity, respectively, such that $[m_{1},m_{2}]=M.$
\end{lemma}

\begin{lemma}
\label{lem2-3}Let $f\in \mathcal{O}(\mathbb{C}^{2},0,0)$\ and let
\begin{equation*}
\mathfrak{M}_{f}=\{m\in \mathbb{N};\;\mathrm{the\;linear\;part\;of\;}f\;%
\mathrm{at}\;0\;\mathrm{has\;periodic\;points\;of\;period\;}m\}.
\end{equation*}%
Then,

\textrm{(i). }For each $m\in \mathbb{N}\backslash \mathfrak{M}_{f}$ such
that $0$ is an isolated fixed point of $f^{m},$%
\begin{equation*}
P_{m}(f,0)=0;
\end{equation*}

\textrm{(ii). }For each positive integer $M$ such that $0$ is an isolated
fixed point of $f^{M},$%
\begin{equation*}
\mu _{f^{M}}(0)=\sum_{\substack{ m\in \mathfrak{M}_{f}  \\ m|M}}P_{m}(f,0).
\end{equation*}
\end{lemma}

\begin{proof}
(i) and (ii) are essentially proved in [\ref{CMY}] (see [\ref{Zh2}] for a
simple proof).
\end{proof}

\begin{remark}
In the previous Lemma, the set $\mathfrak{M}_{f}$ contains at most four
numbers, by Lemma \ref{lem2-5}.
\end{remark}

\begin{lemma}
\label{lem2-4}Let $k$ be a positive integer, let $f$ and $h$ be germs in $%
\mathcal{O}(\mathbb{C}^{2},0,0)$ such that $0$ is an isolated fixed point of
both $f$ and $f^{k}$ and $\det Dh(0)\neq 0,$ and let $g=h\circ f\circ
h^{-1}. $ Then $0$ is still an isolated fixed point of both $g$ and $g^{k}$,
and the following three equalities hold:
\begin{equation*}
\mu _{f^{k}}(0)=\mu _{g^{k}}(0),
\end{equation*}%
\begin{equation*}
P_{k}(f,0)=P_{k}(g,0),
\end{equation*}%
\begin{equation*}
\mathcal{O}_{k}(f,0)=\mathcal{O}_{k}(g,0).
\end{equation*}
\end{lemma}

\begin{proof}
The first equality is well known. The second follows from the first equality
and the definition of Dold's indices. The last then follows from the second
and Corollary \ref{do1}.
\end{proof}

The following result is due to M. Shub and D. Sullivan [\ref{SS}]. It is
also proved in [\ref{Zh1}].

\begin{lemma}
\label{SS-1}Let $m>1$ be a positive integer and let $f\in \mathcal{O}(%
\mathbb{C}^{2},0,0).$ Assume that the origin is an isolated fixed point of $%
f $ and that, for each eigenvalue $\lambda $ of $Df(0),$ either $\lambda =1$
or $\lambda ^{m}\neq 1$. Then the origin is still an isolated fixed point of
$f^{m}$ and%
\begin{equation*}
\mu _{f}(0)=\mu _{f^{m}}(0).
\end{equation*}
\end{lemma}

\begin{lemma}[\protect\cite{LL}]
\label{prod-ind}Let $h_{1}$ and $h_{2}$ be germs in $\mathcal{O}(\mathbb{C}%
^{2},0,0)$. If $0$ is an isolated zero of both $h_{1}$ and $h_{2},$ then the
zero order of $h_{1}\circ h_{2}$ at $0$ equals the product of the zero
orders of $h_{1}$ and $h_{2}$ at $0,$ say, $\pi _{h_{1}\circ h_{2}}(0)=\pi
_{h_{1}}(0)\pi _{h_{2}}(0).$
\end{lemma}

\begin{lemma}[\protect\cite{LL}]
\label{prod1}Let $f=(f_{1},h)$ and $g=(f_{2},h)$ be two germs in $\mathcal{O}%
(\mathbb{C}^{2},0,0).$ If $0$ is an isolated zero of both $f$ and $g,$ then $%
0$ is also an isolated zero of $F=(f_{1}f_{2},h)$ and
\begin{equation*}
\pi _{F}(0)=\pi _{f}(0)+\pi _{h}(0).
\end{equation*}
\end{lemma}

\section{Normal Forms and Iterations of Normal Forms\label{norm}\protect}%
\bigskip

The following lemma is a basic result in the theory of normal forms (see
\cite{AP}, p. 84--85).

\begin{lemma}
\label{nor}Let $f\in \mathcal{O}(\mathbb{C}^{2},0,0)$ and assume that $%
Df(0)=(\lambda _{1},\lambda _{2})$ is a diagonal matrix. Then for any
positive integer $r,$ there exists a polynomial transform
\begin{equation}
(y_{1},y_{2})=H(x_{1},x_{2})=(x_{1},x_{2})+\mathrm{higher\ terms}
\label{tra}
\end{equation}%
of coordinates in a neighborhood of the origin such that each component $%
g_{j}$ of
\begin{equation*}
g=(g_{1},g_{2})=H^{-1}\circ f\circ H
\end{equation*}%
has a power series expansion%
\begin{equation}
g_{j}(x_{1},x_{2})=\lambda
_{j}x_{j}+%
\sum_{_{i_{1}+i_{2}=2}}^{r}c_{i_{1}i_{2}}^{j}x_{1}^{i_{1}}x_{2}^{i_{2}}+%
\mathrm{higher\;terms,\;}j=1,2,  \label{tra-d}
\end{equation}%
in a neighborhood of the origin, where $i_{1}$ and $i_{2}$ are nonnegative
integers and, for $j=1$ and $2,$%
\begin{equation}
c_{i_{1}i_{2}}^{j}\neq 0\ \mathrm{only\ if\ }\lambda _{j}=\lambda
_{1}^{i_{1}}\lambda _{2}^{i_{2}}.  \label{tra-d-1}
\end{equation}
\end{lemma}

\begin{corollary}
\label{tra-d0}Let $H$ be any transform given by the previous lemma. Then for
each $k\in \mathbb{N},$ the $k$ th iteration $%
g^{k}=(g_{1}^{(k)},g_{2}^{(k)}) $ of the germ $g=(g_{1},g_{2})=H^{-1}\circ
f\circ H$ has an expansion similar to (\ref{tra-d}), more precisely, in a
neighborhood of the origin,%
\begin{equation}
g_{j}^{(k)}(x_{1},x_{2})=\lambda
_{j}^{k}x_{j}+%
\sum_{_{i_{1}+i_{2}=2}}^{r}C_{i_{1}i_{2}}^{kj}x_{1}^{i_{1}}x_{2}^{i_{2}}+%
\mathrm{higher\;terms,\;}j=1,2,  \label{tra-d-d}
\end{equation}%
where $i_{1}$ and $i_{2}$ are nonnegative integers and, for $j=1$ and $2,$%
\begin{equation}
C_{i_{1}i_{2}}^{kj}\neq 0\ \mathrm{only\ if\ }\lambda _{j}=\lambda
_{1}^{i_{1}}\lambda _{2}^{i_{2}}.  \label{tra-d-d0}
\end{equation}
\end{corollary}

Before starting the proof, we recall that the notation $o(|z|^{r})$ denotes
any holomorphic function germ whose power series expansion at the origin has
no terms of degrees from $0$ to $r.$ The same notation $o(|z|^{r})$ may
denote different function germs in different places, even in a single
equation.

\begin{proof}
By (\ref{tra-d}), the conclusion holds obviously, except (\ref{tra-d-d0}).
We show that the coefficients $C_{i_{1}i_{2}}^{kj}$ satisfy (\ref{tra-d-d0})
for all $k\in \mathbb{N}$ and $j=1,2.$ This is done by induction on $k.$

Since $g^{1}=(g_{1}^{(1)},g_{2}^{(1)})=(g_{1},g_{2})=g,$ the conclusion
holds for $k=1,$ by the previous lemma. Assume that (\ref{tra-d-d0}) is true
for $k=1,\dots ,l.$ We complete the proof by showing that (\ref{tra-d-d0})
is true for $k=l+1\ $and $j=1,2.$

For $j=1,$ and $k=l+1,$ it is clear by the induction hypothesis that%
\begin{eqnarray*}
&&g_{1}^{(l+1)}(x_{1},x_{2}) \\
&=&g_{1}^{(l)}\circ
g(x_{1},x_{2})=g_{1}^{(l)}(g_{1}(x_{1},x_{2}),g_{2}(x_{1},x_{2})) \\
&=&\lambda
_{1}^{l}g_{1}(x_{1},x_{2})+\sum_{_{i_{1}+i_{2}=2}}^{r}C_{i_{1}i_{2}}^{l1}
\left[ g_{1}(x_{1},x_{2})\right] ^{i_{1}}\left[ g_{2}(x_{1},x_{2})\right]
^{i_{2}}+o(|g(x_{1},x_{2})|^{r}),
\end{eqnarray*}%
and then, writing $o(|g(x_{1},x_{2})|^{r})=o(|x|^{r})$, we have%
\begin{eqnarray}
&&g_{1}^{(l+1)}(x_{1},x_{2})  \label{add} \\
&=&\lambda
_{1}^{l}g_{1}(x_{1},x_{2})+\sum_{_{i_{1}+i_{2}=2}}^{r}C_{i_{1}i_{2}}^{l1}
\left[ g_{1}(x_{1},x_{2})\right] ^{i_{1}}\left[ g_{2}(x_{1},x_{2})\right]
^{i_{2}}+o(|x|^{r}),  \notag
\end{eqnarray}%
and by the induction hypothesis, for each pair $(i_{1},i_{2})$ in the sum in
(\ref{add}),%
\begin{equation}
C_{i_{1}i_{2}}^{l1}\neq 0\;\mathrm{only\ if\ }\lambda _{1}=\lambda
_{1}^{i_{1}}\lambda _{2}^{i_{2}}.  \label{add-1}
\end{equation}

By (\ref{tra-d}), for the first part of the right hand side of (\ref{add}),
we have
\begin{equation*}
\lambda _{1}^{l}g_{1}(x_{1},x_{2})=\lambda _{1}^{l+1}x_{1}+\lambda
_{1}^{l}%
\sum_{_{i_{1}+i_{2}=2}}^{r}c_{i_{1}i_{2}}^{1}x_{1}^{i_{1}}x_{2}^{i_{2}}+%
\mathrm{higher\;terms,}
\end{equation*}%
where each $c_{i_{1}i_{2}}^{1}$ satisfies (\ref{tra-d-1}) for $j=1,$ say, $%
\lambda _{1}^{l}g_{1}(x_{1},x_{2})$ is already in the form of the right hand
side of (\ref{tra-d-d}), together with condition (\ref{tra-d-d0}). Thus, by (%
\ref{add}), to complete the induction for $j=1$, it suffices to show that we
can write the sum%
\begin{equation}
\sigma (x_{1},x_{2})=\sum_{_{i_{1}+i_{2}=2}}^{r}C_{i_{1}i_{2}}^{l1}\left[
g_{1}(x_{1},x_{2})\right] ^{i_{1}}\left[ g_{2}(x_{1},x_{2})\right] ^{i_{2}}
\label{4.7}
\end{equation}%
in (\ref{add}) to be%
\begin{equation}
\sigma (x_{1},x_{2})=\sum_{_{j_{1}+j_{2}=2}}^{r}\mathcal{C}%
_{j_{1}j_{2}}^{l+1,1}x_{1}^{j_{1}}x_{2}^{j_{2}}+o(|x|^{r})\mathrm{,}
\label{h5}
\end{equation}%
such that for each pair $(j_{1},j_{2})$
\begin{equation}
\mathcal{C}_{j_{1}j_{2}}^{l+1,1}\neq 0\mathrm{\ only\ if\ }\ \lambda
_{1}=\lambda _{1}^{j_{1}}\lambda _{2}^{j_{2}}.  \label{h6}
\end{equation}

By (\ref{4.7}) and the expressions of $g_{1}$ and $g_{2}$ in (\ref{tra-d}),
we have
\begin{eqnarray*}
&&\sigma (x_{1},x_{2}) \\
&=&\sum_{_{i_{1}+i_{2}=2}}^{r}C_{i_{1}i_{2}}^{l1}\left[ \lambda
_{1}x_{1}+%
\sum_{s_{1}+s_{2}=2}^{r}c_{s_{1}s_{2}}^{1}x_{1}^{s_{1}}x_{2}^{s_{2}}+o(|x|^{r})%
\right] ^{i_{1}} \\
&&\times \left[ \lambda
_{2}x_{2}+%
\sum_{t_{1}+t_{2}=2}^{r}c_{t_{1}t_{2}}^{2}x_{1}^{t_{1}}x_{2}^{t_{2}}+o(|x|^{r})%
\right] ^{i_{2}} \\
&=&\sum_{_{i_{1}+i_{2}=2}}^{r}C_{i_{1}i_{2}}^{l1}\left[ \lambda
_{1}x_{1}+%
\sum_{s_{1}+s_{2}=2}^{r}c_{s_{1}s_{2}}^{1}x_{1}^{s_{1}}x_{2}^{s_{2}}\right]
^{i_{1}} \\
&&\times \left[ \lambda
_{2}x_{2}+%
\sum_{t_{1}+t_{2}=2}^{r}c_{t_{1}t_{2}}^{2}x_{1}^{t_{1}}x_{2}^{t_{2}}\right]
^{i_{2}}+o(|x|^{r})
\end{eqnarray*}

Thus, $\sigma (x_{1},x_{2})$ is a power series at the origin which is the
sum of terms that is either of the form $cx_{1}^{i_{1}}x_{2}^{i_{2}}$ with $%
i_{1}+i_{2}>r,$ or of the form
\begin{eqnarray}
D_{j_{1}j_{2}}x_{1}^{j_{1}}x_{2}^{j_{2}} &=&C_{i_{1}i_{2}}^{l1}\left[ \left(
\lambda _{1}x_{1}\right) ^{l_{1}}\prod_{(s_{1},s_{2})\in E_{1}}\left(
c_{s_{1}s_{2}}^{1}x_{1}^{s_{1}}x_{2}^{s_{2}}\right) ^{l_{s_{1}s_{2}}}\right]
\label{h} \\
&&\times \left[ \left( \lambda _{2}x_{2}\right) ^{l_{2}^{\prime
}}\prod_{(t_{1},t_{2})\in E_{2}}\left(
c_{t_{1}t_{2}}^{2}x_{1}^{t_{1}}x_{2}^{t_{2}}\right) ^{l_{t_{1}t_{2}}^{\prime
}}\right] ,  \notag
\end{eqnarray}%
where $l_{1}$, $l_{2}^{\prime },$ $l_{s_{1}s_{2}}$ and $l_{t_{1}t_{2}}^{%
\prime }$ are nonnegative integers, and $E_{1}$ and $E_{2}$ are sets of some
pairs of nonnegative integers such that
\begin{equation}
l_{1}+\sum_{(s_{1},s_{2})\in E_{1}}l_{s_{1}s_{2}}=i_{1},\ \mathrm{and}\
l_{2}^{\prime }+\sum_{(t_{1},t_{2})\in E_{2}}l_{t_{1}t_{2}}^{\prime }=i_{2}
\label{h1}
\end{equation}%
\begin{equation}
D_{j_{1}j_{2}}=\lambda _{1}^{l_{1}}\lambda _{2}^{l_{2}^{\prime
}}C_{i_{1}i_{2}}^{l1}\left[ \prod_{(s_{1},s_{2})\in E_{1}}\left(
c_{s_{1}s_{2}}^{1}\right) ^{l_{s_{1}s_{2}}}\right] \left[
\prod_{(t_{1},t_{2})\in E_{2}}\left( c_{t_{1}t_{2}}^{2}\right)
^{l_{t_{1}t_{2}}^{\prime }}\right] ,  \label{h2}
\end{equation}%
\begin{equation}
j_{1}=l_{1}+\sum_{(s_{1},s_{2})\in
E_{1}}s_{1}l_{s_{1}s_{2}}+\sum_{(t_{1},t_{2})\in
E_{2}}t_{1}l_{t_{1}t_{2}}^{\prime },  \label{h3}
\end{equation}%
and%
\begin{equation}
j_{2}=l_{2}^{\prime }+\sum_{(s_{1},s_{2})\in
E_{1}}s_{2}l_{s_{1}s_{2}}+\sum_{(t_{1},t_{2})\in
E_{2}}t_{2}l_{t_{1}t_{2}}^{\prime }.  \label{h4}
\end{equation}%
We show that
\begin{equation*}
D_{j_{1}j_{2}}\neq 0\mathrm{\ only\ if\ }\lambda _{1}=\lambda
_{1}^{j_{1}}\lambda _{2}^{j_{2}}.
\end{equation*}

By (\ref{h2}), (\ref{tra-d-1}) and (\ref{add-1}), it is clear that $%
D_{j_{1}j_{2}}\neq 0$ only if

\begin{equation*}
\lambda _{1}=\lambda _{1}^{i_{1}}\lambda _{2}^{i_{2}},
\end{equation*}%
\begin{equation*}
\lambda _{1}=\lambda _{1}^{s_{1}}\lambda _{2}^{s_{2}}\;\mathrm{for\ all\ }%
(s_{1},s_{2})\in E_{1},
\end{equation*}%
and%
\begin{equation*}
\lambda _{2}=\lambda _{1}^{t_{1}}\lambda _{2}^{t_{2}}\;\mathrm{for\ all\ }%
(t_{1},t_{2})\in E_{2}.
\end{equation*}%
Therefore, if $D_{j_{1}j_{2}}\neq 0,$ then we have, by (\ref{h1})--(\ref{h4}%
), that
\begin{eqnarray*}
\lambda _{1}^{j_{1}}\lambda _{2}^{j_{2}} &=&\lambda
_{1}^{l_{1}+\sum_{(s_{1},s_{2})\in
E_{1}}s_{1}l_{s_{1}s_{2}}+\sum_{(t_{1},t_{2})\in
E_{2}}t_{1}l_{t_{1}t_{2}}^{\prime }} \\
&&\times \lambda _{2}^{l_{2}^{\prime }+\sum_{(s_{1},s_{2})\in
E_{1}}s_{2}l_{s_{1}s_{2}}+\sum_{(t_{1},t_{2})\in
E_{2}}t_{2}l_{t_{1}t_{2}}^{\prime }} \\
&=&\left[ \lambda _{1}^{l_{1}}\prod_{(s_{1},s_{2})\in E_{1}}\left( \lambda
_{1}^{s_{1}}\lambda _{2}^{s_{2}}\right) ^{l_{s_{1}s_{2}}}\right] \left[
\lambda _{2}^{l_{2}^{\prime }}\prod_{(t_{1},t_{2})\in E_{2}}\left( \lambda
_{1}^{t_{1}}\lambda _{2}^{t_{2}}\right) ^{l_{t_{1}t_{2}}^{\prime }}\right] \\
&=&\left[ \lambda _{1}^{l_{1}}\prod_{(s_{1},s_{2})\in E_{1}}\left( \lambda
_{1}\right) ^{l_{s_{1}s_{2}}}\right] \left[ \lambda _{2}^{l_{2}^{\prime
}}\prod_{(t_{1},t_{2})\in E_{2}}(\lambda _{2})^{l_{t_{1}t_{2}}^{\prime }}%
\right] \\
&=&\left[ \lambda _{1}^{l_{1}+\sum_{(s_{1},s_{2})\in E_{1}}l_{s_{1}s_{2}}}%
\right] \left[ \lambda _{2}^{l_{2}^{\prime }+\sum_{(t_{1},t_{2})\in
E_{2}}l_{t_{1}t_{2}}^{\prime }}\right] \\
&=&\lambda _{1}^{i_{1}}\lambda _{2}^{i_{2}}=\lambda _{1}.
\end{eqnarray*}%
In other words, $D_{j_{1}j_{2}}\neq 0$ only if $\lambda _{1}=\lambda
_{1}^{j_{1}}\lambda _{2}^{j_{2}}.$ Thus, we can write $\sigma (x_{1},x_{2})$
to be (\ref{h5}) with (\ref{h6}), and then, we have proved that (\ref%
{tra-d-d0}) is true for $k=l+1$ and $j=1.$ For the same reason, (\ref%
{tra-d-d0}) is true for $k=l+1$ and $j=2.$ The induction is complete.
\end{proof}

\section{Computing Zero Orders Via Cronin's Theorem\label{cronin-1}}

Since fixed point indices are defined via zero orders, to prove the main
theorem, it is useful to compute the zero orders of some special germs of
holomorphic mappings. The following Cronin's theorem plays an important role
in our computations.

\begin{theorem}[\textbf{Cronin \protect\cite{Cro}}]
Let $f\in \mathcal{O}(\mathbb{C}^{2},0,0)$ be given by%
\begin{equation*}
f(x_{1},x_{2})=\left(
P_{m_{1}}(x_{1},x_{2})+o(|x|^{m_{1}}),Q_{m_{2}}(x_{1},x_{2})+o(|x|^{m_{2}}%
\right) ),
\end{equation*}%
where $x=(x_{1},x_{2}),$ $P_{m_{1}}$ and $Q_{m_{2}}$\ are homogeneous
polynomials of degrees $m_{1}$ and $m_{2},$ respectively,\ in $x_{1}$ and $%
x_{2}$. If the origin $0$ is an isolated solution of the system
\begin{equation}
\left\{
\begin{array}{c}
P_{m_{1}}(x_{1},x_{2})=0, \\
Q_{m_{2}}(x_{1},x_{2})=0,%
\end{array}%
\right.  \label{cr}
\end{equation}%
then $0$\ is an isolated zero of the germ $f$ with zero order
\begin{equation*}
\pi _{f}(0)=m_{1}m_{2}.
\end{equation*}%
If $0$ is an isolated zero of $f$ but is not an isolated solution of the
system (\ref{cr}), then%
\begin{equation*}
\pi _{f}(0)>m_{1}m_{2}.
\end{equation*}
\end{theorem}

We now apply Cronin's theorem to some special cases.

\begin{corollary}
\label{c1}Let $g=(g_{1},g_{2})\in \mathcal{O}(\mathbb{C}^{2},0,0)$ be given
by%
\begin{equation*}
\left\{
\begin{array}{l}
g_{1}(x_{1},x_{2})=x_{1}^{m_{1}+1}(a_{1}+o(1))+x_{2}^{d}(a_{2}+o(1)), \\
g_{2}(x_{1},x_{2})=x_{1}^{m_{1}}x_{2}(a_{3}+o(1))+x_{2}^{d+1}O(1),%
\end{array}%
\right.
\end{equation*}%
where $d$ and $m_{1}$ are positive integers, and $a_{i}$ are constants.

If the origin is an isolated zero of $g,$ then%
\begin{equation*}
\pi _{g}(0)\geq dm_{1}+m_{1}+1.
\end{equation*}

If $a_{1}\neq 0,a_{2}\neq 0$ and $a_{3}\neq 0,$ then the origin is an
isolated zero of $g$ with%
\begin{equation*}
\pi _{g}(0)=dm_{1}+m_{1}+1.
\end{equation*}
\end{corollary}

\begin{proof}
Let $h\in \mathcal{O}(\mathbb{C}^{2},0,0)$ be given by
\begin{equation*}
h(y_{1},y_{2})=(y_{1}^{d},y_{2}^{m_{1}+1}).
\end{equation*}%
then the two component of $g\circ h=(g_{1}\circ h,g_{2}\circ h)$ have the
expressions%
\begin{equation}
\left\{
\begin{array}{l}
g_{1}\circ h(y_{1},y_{2})=a_{1}y_{1}^{d(m_{1}+1)}+a_{2}y_{2}^{d(m_{1}+1)}+%
\mathrm{higher\ terms,} \\
g_{2}\circ h(y_{1},y_{2})=a_{3}y_{1}^{dm_{1}}y_{2}^{m_{1}+1}+\mathrm{higher\
terms}.%
\end{array}%
\right.  \label{j14}
\end{equation}%
If the origin is an isolated zero of $g,$ then it is also an isolated zero
of $g\circ h,$ and then, by Cronin's theorem,
\begin{equation}
\pi _{g\circ h}(0)\geq d(m_{1}+1)(dm_{1}+m_{1}+1),  \label{h7}
\end{equation}%
and, considering that $\pi _{h}(0)=d(m_{1}+1),$ by Lemma \ref{prod-ind} we
have
\begin{equation}
\pi _{g}(0)=\pi _{g\circ h}(0)/\pi _{h}(0)\geq dm_{1}+m_{1}+1.  \label{h8}
\end{equation}

On the other hand, it is clear that when $a_{1},a_{2}$ and $a_{3}$ are all
nonzero, $0$ is an isolated solution of the system
\begin{equation*}
\left\{
\begin{array}{c}
a_{1}y_{1}^{d(m_{1}+1)}+a_{2}y_{2}^{d(m_{1}+1)}=0, \\
a_{3}y_{1}^{dm_{1}}y_{2}^{m_{1}+1}=0,%
\end{array}%
\right.
\end{equation*}%
and then by (\ref{j14}) and Cronin's theorem, $0$ is an isolated zero of $%
g\circ h$ and the equality in (\ref{h7}) holds, and then $0$ is also an
isolated zero of $g$ and the equality in (\ref{h8}) holds again\label{?1}.
This completes the proof.
\end{proof}

\begin{corollary}
\label{c2}Let $g=(g_{1},g_{2})\in \mathcal{O}(\mathbb{C}^{2},0,0)$ be given
by%
\begin{equation*}
\left\{
\begin{array}{l}
g_{1}(x_{1},x_{2})=x_{1}^{m}(a+o(1))+x_{2}O(1), \\
g_{2}(x_{1},x_{2})=x_{2}(b+o(1))+x_{1}^{m}o(1),%
\end{array}%
\right.
\end{equation*}%
where $a$ and $b$ are constants.

If the origin is an isolated zero of $g,$ then%
\begin{equation*}
\pi _{g}(0)\geq m.
\end{equation*}

If $a\neq 0$ and $b\neq 0$, then the origin is an isolated zero of $g$ with
\begin{equation*}
\pi _{g}(0)=m.
\end{equation*}
\end{corollary}

\begin{proof}
Let $h\in \mathcal{O}(\mathbb{C}^{2},0,0)$ be given by
\begin{equation*}
h(y_{1},y_{2})=(y_{1},y_{2}^{m}).
\end{equation*}%
Then the germ $g\circ h=(g_{1}\circ h,g_{2}\circ h)$ has the expression

\begin{equation}
\left\{
\begin{array}{l}
g_{1}\circ h(y_{1},y_{2})=ay_{1}^{m}+a_{0}y_{2}^{m}+\mathrm{higher\ terms,}
\\
g_{2}\circ h(y_{1},y_{2})=by_{2}^{m}+\mathrm{higher\ terms},%
\end{array}%
\right.  \label{j4}
\end{equation}%
for some constant $a_{0}.$

If the origin is an isolated zero of $g,$ then it is also an isolated zero
of the germ $g\circ h$, and by (\ref{j4}) and Cronin's theorem,
\begin{equation*}
\pi _{g\circ h}(0)\geq m^{2},
\end{equation*}%
and the equality holds if $a\neq 0$ and $b\neq 0;$ and then by Lemma \ref%
{prod-ind} and by the fact $\pi _{h}(0)=m,$%
\begin{equation*}
\pi _{g}(0)=\pi _{g\circ h}(0)/\pi _{h}(0)\geq m,
\end{equation*}%
and the equality holds if $a\neq 0$ and $b\neq 0.$ On the other hand, by (%
\ref{j4}) and Cronin's theorem, when $a\neq 0$ and $b\neq 0,$ $0$ is an
isolated zero of $g\circ h,$ and then $0$ is also an isolated zero of $g.$
This completes the proof.
\end{proof}

\begin{corollary}
\label{c3}Let $f=(f_{1},f_{2})\in \mathcal{O}(\mathbb{C}^{2},0,0)$ be given
by%
\begin{eqnarray*}
f_{1}(x_{1},x_{2}) &=&x_{1}^{n_{1}}O(1)+x_{2}^{n_{2}}O(1), \\
f_{2}(x_{1},x_{2}) &=&x_{1}^{n_{1}}O(1)+x_{2}^{n_{2}}O(1),
\end{eqnarray*}%
where $n_{1}$ and $n_{2}$ are positive integers. Assume that $0$ is an
isolated zero of $f.$ Then
\begin{equation*}
\pi _{f}(0)\geq n_{1}n_{2}.
\end{equation*}
\end{corollary}

\begin{proof}
Consider $h\in \mathcal{O}(\mathbb{C}^{2},0,0)$ given by $%
h(y_{1},y_{2})=(y_{1}^{n_{2}},y_{2}^{n_{1}}).$ Then $0$ is an isolated zero
of the germ $f\circ h=(f_{1}\circ h,f_{2}\circ h)$ and
\begin{eqnarray*}
f_{1}\circ h(y_{1},y_{2}) &=&y_{1}^{n_{1}n_{2}}O(1)+y_{2}^{n_{1}n_{2}}O(1),
\\
f_{2}\circ h(y_{1},y_{2}) &=&y_{1}^{n_{1}n_{2}}O(1)+y_{2}^{n_{1}n_{2}}O(1).
\end{eqnarray*}%
Thus, by Cronin's theorem we have $\pi _{f\circ h}(0)\geq
n_{1}^{2}n_{2}^{2}, $ and then by the fact $\pi _{h}(0)=n_{1}n_{2}$ and by
Lemma \ref{prod-ind}, we have $\pi _{f}(0)=\pi _{f\circ h}(0)/\pi
_{h}(0)\geq n_{1}n_{2}.$
\end{proof}

\begin{corollary}
\label{c4}Let $f=(f_{1},f_{2})\in \mathcal{O}(\mathbb{C}^{2},0,0)$ be given
by%
\begin{eqnarray*}
f_{1}(x_{1},x_{2}) &=&x_{1}^{m_{1}}, \\
f_{2}(x_{1},x_{2}) &=&x_{1}^{m_{1}}O(1)+x_{2}^{rm_{2}}O(1),
\end{eqnarray*}%
where $r,m_{1}$ and $m_{2}$ are positive integers. Assume that $0$ is an
isolated zero of $f.$ Then
\begin{equation*}
\pi _{f}(0)\geq rm_{1}m_{2}.
\end{equation*}
\end{corollary}

\begin{proof}
It follows from the previous corollary, by taking $n_{1}=m_{1}$ and $%
n_{2}=rm_{2}.$
\end{proof}

\begin{corollary}
\label{c5}Let $f=(f_{1},f_{2})\in \mathcal{O}(\mathbb{C}^{2},0,0)$ be given
by%
\begin{eqnarray*}
f_{1}(x_{1},x_{2}) &=&x_{1}^{n_{1}}\left[
x_{1}^{rn_{1}}O(1)+x_{2}^{n_{2}}O(1)\right] , \\
f_{2}(x_{1},x_{2}) &=&x_{2}^{n_{2}}\left[
x_{1}^{n_{1}}O(1)+x_{2}^{rn_{2}}O(1)\right] ,
\end{eqnarray*}%
where $r,$ $n_{1}$ and $n_{2}$ are positive integers. Assume that $0$ is an
isolated zero of $f.$ Then%
\begin{equation*}
\pi _{f}(0)\geq 2n_{1}n_{2}+2rn_{1}n_{2}.
\end{equation*}
\end{corollary}

\begin{proof}
By Lemma \ref{prod1}, the zero order $\pi _{f}(0)$ equals the sum of the
zero orders of the four germs in $\mathcal{O}(\mathbb{C}^{2},0,0)$ given by%
\begin{eqnarray*}
g_{1}(x_{1},x_{2}) &=&(x_{1}^{n_{1}},x_{2}^{n_{2}}), \\
g_{2}(x_{1},x_{2}) &=&(x_{1}^{n_{1}},x_{1}^{n_{1}}O(1)+x_{2}^{rn_{2}}O(1)),
\\
g_{3}(x_{1},x_{2}) &=&(x_{1}^{rn_{1}}O(1)+x_{2}^{n_{2}}O(1),x_{2}^{n_{2}}),
\\
g_{4}(x_{1},x_{2})
&=&(x_{1}^{rn_{1}}O(1)+x_{2}^{n_{2}}O(1),x_{1}^{n_{1}}O(1)+x_{2}^{rn_{2}}O(1)).
\end{eqnarray*}%
By Cronin's theorem we have $\pi _{g_{1}}(0)\geq n_{1}n_{2};$ by Corollary %
\ref{c4} we have $\pi _{g_{2}}(0)\geq rn_{1}n_{2}$ and $\pi _{g_{3}}(0)\geq
rn_{1}n_{2};$ and by Corollary \ref{c3} we have $\pi _{g_{4}}(0)\geq
n_{1}n_{2}.$ Thus, we have
\begin{equation*}
\pi _{f}(0)\geq 2n_{1}n_{2}+2rn_{1}n_{2}.
\end{equation*}
\end{proof}

\begin{lemma}
\label{ad-2-1}Let $f=(f_{1},f_{2})$ and $g=(g_{1},g_{2})$ be germs in $%
\mathcal{O}(\mathbb{C}^{2},0,0)$ and assume that $A=(a_{ij})$ is a $2\times
2 $ matrix whose elements $a_{ij}$ are germs of holomorphic functions at the
origin of $\mathbb{C}^{2},$ with $\det A(0)\neq 0$. If%
\begin{equation*}
(f_{1},f_{2})=(g_{1},g_{2})A=(g_{1}a_{11}+g_{2}a_{21},g_{1}a_{12}+g_{2}a_{22}),
\end{equation*}%
and the origin is an isolated zero of $g,$ then the origin is also an
isolated zero of $f$ and
\begin{equation*}
\pi _{f}(0)=\pi _{g}(0).
\end{equation*}
\end{lemma}

\begin{proof}
By the assumption, there exists a ball $B$ centered at the origin in $%
\mathbb{C}^{2}$ such that the origin is the unique zero of $g$ in $\overline{%
B},$ $A$ is well defined on $\overline{B}$ and
\begin{equation}
\det A(x_{1},x_{2})\neq 0,\ \ (x_{1},x_{2})\in B.  \label{h9}
\end{equation}

Then there exists a regular value $\varepsilon =(\varepsilon
_{1},\varepsilon _{2})$ of $g,$ which can be chosen close to the origin
arbitrarily, such that $g^{-1}(\varepsilon )\cap B$ contains exactly $\pi
_{g}(0)$ distinct points $(a_{1},b_{1}),\dots ,(a_{\pi _{g}(0)},b_{\pi
_{g}(0)})$ in $B.$ Thus, we have%
\begin{equation*}
f(a_{i},b_{i})-(\varepsilon _{1},\varepsilon _{2})A(a_{i},b_{i})=0,i=1,\dots
,\pi _{g}(0).
\end{equation*}%
In other words, $f_{\varepsilon }(x_{1},x_{2})=f(x_{1},x_{2})-(\varepsilon
_{1},\varepsilon _{2})A(x_{1},x_{2})$ has $\pi _{g}(0)$ distinct zeros in $%
B. $ It is clear that $f_{\varepsilon }(x_{1},x_{2})$ is a germ in $\mathcal{%
O}(\mathbb{C}^{2},0,0)$ converging to $f$ uniformly on $\overline{B}$ as $%
\varepsilon =(\varepsilon _{1},\varepsilon _{2})\rightarrow 0,$ and on the
other hand, the origin is also the unique zero of $f$ in $\overline{B}$ as
well. Thus, by Rouche's theorem, we have%
\begin{equation*}
\pi _{f}(0)\geq \pi _{g}(0).
\end{equation*}%
But by (\ref{h9}), the inverse $A^{-1}$ of the matrix $A$ is well defined on
$\overline{B},$ which is again a matrix of holomorphic functions, and $%
g=fA^{-1}.$ Thus, for the same reason we have $\pi _{f}(0)\leq \pi _{g}(0).$
This completes the proof.
\end{proof}

\begin{corollary}
\label{ad-1}Let $f=(f_{1},f_{2})\in \mathcal{O}(\mathbb{C}^{2},0,0)$, let $h$
be a holomorphic function germ at the origin, and let $g\in \mathcal{O}(%
\mathbb{C}^{2},0,0)$ be given by%
\begin{equation*}
g=(f_{1},f_{2}+hf_{1})=(f_{1},f_{2})\left(
\begin{array}{cc}
1 & \ h \\
0 & \ 1%
\end{array}%
\right) .
\end{equation*}%
If the origin is an isolated zero of $f,$ then it is an isolated zero of $g$
and
\begin{equation*}
\pi _{g}(0)=\pi _{f}(0).
\end{equation*}
\end{corollary}

\begin{lemma}
\label{dd4}Let $d>1,n_{1}>1$ and $n_{2}>1$ be positive integers and let $%
f=(f_{1},f_{2})\in \mathcal{O}(\mathbb{C}^{2},0,0)$ be given by

\begin{equation}
\left\{
\begin{array}{rl}
f_{1}(x_{1},x_{2})= &
x_{1}^{dn_{1}+1}a_{11}(x_{1},x_{2})+x_{1}^{n_{1}+1}x_{2}^{n_{2}}O(1) \\
& +x_{1}x_{2}^{dn_{2}}O(1)+x_{2}^{2dn_{2}+1}a_{12}(x_{1},x_{2}),\medskip \\
f_{2}(x_{1},x_{2})= &
x_{1}^{dn_{1}}x_{2}O(1)+x_{1}^{n_{1}}x_{2}^{n_{2}+1}O(1) \\
& +x_{2}^{dn_{2}+1}a_{22}(x_{1},x_{2})+x_{1}^{2dn_{1}+1}a_{21}(x_{1},x_{2}),%
\end{array}%
\right.  \label{dd}
\end{equation}%
where $a_{ij}=a_{ij}(x_{1},x_{2})$ are holomorphic function germs at the
origin. Assume that the origin is an isolated zero of $f.$ Then
\begin{equation}
\pi _{f}(0)\geq 2dn_{1}n_{2}+dn_{1}+dn_{2}+1.  \label{dd1}
\end{equation}
\end{lemma}

\begin{proof}
We first assume%
\begin{equation}
a_{11}(0)\neq 0,a_{22}(0)\neq 0.  \label{dd3}
\end{equation}%
Then $[a_{11}(x_{1},x_{2})]^{-1}$ and $[a_{22}(x_{1},x_{2})]^{-1}$ are also
holomorphic function germs at the origin, and we can reduce the germ $f$
into a simpler germ $h=(h_{1},h_{2})$ of the form%
\begin{equation}
\left\{
\begin{array}{c}
h_{1}(x_{1},x_{2})=x_{1}r_{1}(x_{1},x_{2}), \\
h_{2}(x_{1},x_{2})=x_{2}r_{2}(x_{1},x_{2}).%
\end{array}%
\right.  \label{ddd1}
\end{equation}%
with
\begin{equation}
\left\{
\begin{array}{c}
r_{1}(x_{1},x_{2})=x_{1}^{dn_{1}}a_{11}(x_{1},x_{2})+x_{1}^{n_{1}}x_{2}^{n_{2}}O(1)+x_{2}^{dn_{2}}O(1),
\\
r_{2}(x_{1},x_{2})=x_{1}^{dn_{1}}O(1)+x_{1}^{n_{1}}x_{2}^{n_{2}}O(1)+x_{2}^{dn_{2}}a_{22}(x_{1},x_{2}),%
\end{array}%
\right.  \label{j9}
\end{equation}%
and with
\begin{equation}
\pi _{f}(0)=\pi _{h}(0).  \label{j10}
\end{equation}

By Corollary \ref{ad-1}, $\pi _{f}(0)$ equals $\pi _{g}(0)$ for the germ $%
g=(g_{1},g_{2})\in \mathcal{O}(\mathbb{C}^{2},0,0)$ that is given by
\begin{eqnarray*}
g_{1}(x_{1},x_{2}) &=&f_{1}(x_{1},x_{2}), \\
g_{2}(x_{1},x_{2})
&=&f_{2}(x_{1},x_{2})-x_{1}^{dn_{1}}[a_{11}(x_{1},x_{2})]^{-1}a_{21}(x_{1},x_{2})f_{1}(x_{1},x_{2}),
\end{eqnarray*}%
Since $n_{1}>1,n_{2}>1$ and $d>1,$ by (\ref{dd}) we can write%
\begin{equation*}
x_{1}^{dn_{1}}[a_{11}(x_{1},x_{2})]^{-1}a_{21}(x_{1},x_{2})f_{1}(x_{1},x_{2})=x_{1}^{2dn_{1}+1}a_{21}(x_{1},x_{2})+x_{1}^{dn_{1}}x_{2}O(1),
\end{equation*}%
and then again by (\ref{dd}) we have%
\begin{equation}
g_{2}(x_{1},x_{2})=x_{1}^{dn_{1}}x_{2}O(1)+x_{1}^{n_{1}}x_{2}^{n_{2}+1}O(1)+x_{2}^{dn_{2}+1}a_{22}(x_{1},x_{2}).
\label{j13}
\end{equation}%
Again by Corollary \ref{ad-1}, $\pi _{g}(0)$ equals $\pi _{h}(0)$ for the
germ $h=(h_{1},h_{2})$ given by%
\begin{eqnarray*}
h_{1}(x_{1},x_{2})
&=&g_{1}(x_{1},x_{2})-x_{2}^{dn_{2}}[a_{22}(x_{1},x_{2})]^{-1}a_{12}(x_{1},x_{2})g_{2}(x_{1},x_{2}),
\\
h_{2}(x_{1},x_{2}) &=&g_{2}(x_{1},x_{2}),
\end{eqnarray*}%
and, by (\ref{j13}) and the expression of $g_{1}=f_{1}$ in (\ref{dd}), it is
easy to see that $h_{1}(x_{1},x_{2})$ has the expression%
\begin{equation*}
h_{1}(x_{1},x_{2})=x_{1}^{dn_{1}+1}a_{11}(x_{1},x_{2})+x_{1}^{n_{1}+1}x_{2}^{n_{2}}O(1)+x_{1}x_{2}^{dn_{2}}O(1),
\end{equation*}%
and then the germ $h=(h_{1},h_{2})=(h_{1},g_{2})$ has the expression (\ref%
{ddd1}), such that (\ref{j9}) and (\ref{j10}) hold.

By (\ref{dd3}), repeating the above arguments, we can reduce the germ $%
r=(r_{1},r_{2})\in \mathcal{O}(\mathbb{C}^{2},0,0)$ into a further simpler
germ $s=(s_{1},s_{2})$ with the expression%
\begin{eqnarray*}
s_{1}(x_{1},x_{2})
&=&x_{1}^{dn_{1}}O(1)+x_{1}^{n_{1}}x_{2}^{n_{2}}O(1)=x_{1}^{n_{1}}\left[
x_{1}^{(d-1)n_{1}}O(1)+x_{2}^{n_{2}}O(1)\right] , \\
s_{2}(x_{1},x_{2})
&=&x_{1}^{n_{1}}x_{2}^{n_{2}}O(1)+x_{2}^{dn_{2}}O(1)=x_{2}^{n_{2}}\left[
x_{1}^{n_{1}}O(1)+x_{2}^{(d-1)n_{2}}O(1)\right] ,
\end{eqnarray*}%
such that $\pi _{r}(0)=\pi _{s}(0).$ By Corollary \ref{c5}, we have%
\begin{equation*}
\pi _{s}(0)\geq 2n_{1}n_{2}+2(d-1)n_{1}n_{2}=2dn_{1}n_{2},
\end{equation*}%
and then%
\begin{equation*}
\pi _{r}(0)=\pi _{s}(0)\geq 2dn_{1}n_{2}.
\end{equation*}

On the other hand, by Corollary \ref{c2}, for the germs in $\mathcal{O}(%
\mathbb{C}^{2},0,0)$ given by%
\begin{eqnarray*}
k_{1}(x_{1},x_{2}) &=&(x_{1},x_{2}), \\
k_{2}(x_{1},x_{2}) &=&(r_{1}(x_{1},x_{2}),x_{2}), \\
k_{3}(x_{1},x_{2}) &=&(x_{1},r_{2}(x_{1},x_{2})),
\end{eqnarray*}%
we have%
\begin{equation*}
\pi _{k_{1}}(0)=1,\pi _{k_{2}}(0)=dn_{1},\pi _{k_{3}}(0)=dn_{2}.
\end{equation*}%
Thus, by Lemma \ref{prod1} and (\ref{ddd1}),%
\begin{eqnarray*}
\pi _{h}(0) &=&\pi _{k_{1}}(0)+\pi _{k_{2}}(0)+\pi _{k_{3}}(0)+\pi _{r}(0) \\
&\geq &1+dn_{1}+dn_{2}+2dn_{1}n_{2},
\end{eqnarray*}%
which implies (\ref{dd1}), by (\ref{j10}).

If (\ref{dd3}) fails, then we consider the germ $f_{\varepsilon }\in
\mathcal{O}(\mathbb{C}^{2},0,0)$ that is obtain from $f$ by just replacing $%
a_{ii}(x_{1},x_{2})$ with $a_{ii}(x_{1},x_{2})+\varepsilon ,i=1,2.$ Then $%
f_{\varepsilon }$ converges to $f$ uniformly in a neighborhood of the
origin, $f_{\varepsilon }$ has the form of (\ref{dd}) and for sufficiently
small $\varepsilon ,$ $0$ is an isolated zero of $f_{\varepsilon }$ by Rouch%
\'{e}'s Theorem, and $a_{ii}(0,0)+\varepsilon \neq 0$ for $i=1$ and $2.$
Thus, the above arguments are applied to such $f_{\varepsilon }$ if $%
\varepsilon $ is small enough. In other word, for sufficiently small $%
\varepsilon ,$ we have%
\begin{equation*}
\pi _{f_{\varepsilon }}(0)\geq 1+dn_{1}+dn_{2}+2dn_{1}n_{2},
\end{equation*}%
and then (\ref{dd1}) follows from Corollary \ref{ad-2-0}.
\end{proof}

By Cronin's theorem, one can prove the following result.

\begin{proposition}
\label{c8}Let $f=(f_{1},f_{2})\in \mathcal{O}(\mathbb{C}^{2},0,0)$ be given
by%
\begin{eqnarray*}
f_{1}(x_{1},x_{2}) &=&\lambda
_{1}x_{1}+x_{1}(a_{11}x_{1}^{m_{1}}+a_{12}x_{2}^{m_{2}}), \\
f_{2}(x_{1},x_{2}) &=&\lambda
_{2}x_{2}+x_{2}(a_{21}x_{1}^{m_{1}}+a_{22}x_{2}^{m_{2}}),
\end{eqnarray*}%
where $\lambda _{1},\lambda _{2}$ are primitive $m_{1}$-th$\ $and $m_{2}$-th
roots of unity, respectively, $m_{1}$ and $m_{2}$ are positive integers that
are relatively prime.

If $a_{11}\neq 0,$ $a_{22}\neq 0$ and $\det (a_{ij})\neq 0,$ then the origin
is an isolated fixed point of $f^{m_{1}},f^{m_{2}}$ and $f^{m_{1}m_{2}}$ and
the following formulae hold.%
\begin{equation*}
\mu _{f^{m_{1}}}(0)=(m_{1}+1),
\end{equation*}%
\begin{equation*}
\mu _{f^{m_{2}}}(0)=(m_{2}+1),
\end{equation*}%
\begin{equation*}
P_{m_{1}m_{2}}(f,0)=m_{1}m_{2}.
\end{equation*}
\end{proposition}

\begin{proof}
When $m_{1}$ and $m_{2}$ are distinct primes, this is proved in [\ref{Zh2}].
But in general, the proof is exactly the same.
\end{proof}

\section{Proof of the Main Theorem: (B)$\Rightarrow $(A) \label{cronin}}

In this section, we deduce (A) from (B) in the main theorem.

Assume that $M>1$ is an integer and $A$ is a matrix that satisfies (B) and
let $f\ $be a germ in $\mathcal{O}(\mathbb{C}^{2},0,0)$ such that%
\begin{equation*}
Df(0)=A,
\end{equation*}%
and that the origin is an isolated fixed point of $f^{M}.$ We shall show that%
\begin{equation}
\mathcal{O}_{M}(f,0)\geq 2.  \label{tar}
\end{equation}

By Lemma \ref{lem2-4} and the assumption in (B), we may assume that%
\begin{equation*}
A=\left(
\begin{array}{cc}
\lambda _{1} & 0 \\
0 & \lambda _{2}%
\end{array}%
\right) ,
\end{equation*}%
where $\lambda _{1}$ and $\lambda _{2}$ are primitive $m_{1}$ th and $m_{2}$
th roots of unity, respectively, and one of the following conditions holds.

\textrm{(b1) }$m_{1}=m_{2}=M$ and $\lambda _{1}=\lambda _{2}.$

\textrm{(b2) }$m_{1}=m_{2}=M$ and there exists positive integers $1<\alpha
<M $ and $1<\beta <M$ such that%
\begin{equation}
\lambda _{1}^{\beta }=\lambda _{2},\lambda _{2}^{\alpha }=\lambda
_{1},\alpha \beta >M+1.  \label{aid}
\end{equation}

\textrm{(b3) }$m_{1}|m_{2},$ $m_{2}=M$, and $\lambda _{2}^{m_{2}/m_{1}}\neq
\lambda _{1}.$

\textrm{(b4) }$M=[m_{1},m_{2}],(m_{1},m_{2})>1$ and $\max \{m_{1},m_{2}\}<M.$

Then, in any case from (b1) to (b4), the origin is a simple fixed point of $%
f $, and then by Lemma \ref{lem2-a}, we have%
\begin{equation}
\mu _{f}(0)=P_{1}(f,0)=1.  \label{all}
\end{equation}

We show that any one of the four conditions from (b1) to (b4) deduces (\ref%
{tar}), and divide the proof into four parts.\medskip

\begin{description}
\item[\textbf{Part 1}] \textbf{(b1) }$\Rightarrow $\textbf{\ (\ref{tar}%
).\medskip }
\end{description}

\begin{proof}
In case (b1), we may assume $\lambda _{1}=\lambda _{2}=\lambda $ and $%
\lambda $ is a primitive $M$ th root of unity. Then,
\begin{equation*}
Df(0)=\left(
\begin{array}{cc}
\lambda & 0 \\
0 & \lambda%
\end{array}%
\right) .
\end{equation*}

By Lemma \ref{nor}, there exists a polynomial transform $%
(y_{1},y_{2})=H(x_{1},x_{2})$ in the form of (\ref{tra}), such that each
component of $g=(g_{1},g_{2})=H^{-1}\circ f\circ H$ has the expression%
\begin{equation}
g_{j}(x_{1},x_{2})=\lambda
x_{j}+%
\sum_{_{i_{1}+i_{2}=2}}^{M}c_{i_{1}i_{2}}^{j}x_{1}^{i_{1}}x_{2}^{i_{2}}+o(|x|^{M}),j=1,2,
\label{nor1}
\end{equation}%
in a neighborhood of the origin, where the sum extends over all $2$-tuples $%
(i_{1},i_{2})$ of nonnegative integers with%
\begin{equation*}
2\leq i_{1}+i_{2}\leq M+1\ \mathrm{and\ }\lambda =\lambda ^{i_{1}+i_{2}},
\end{equation*}%
which implies that $M|(i_{1}+i_{2}-1),$ since $\lambda $ is a primitive $M$
th root of unity. Thus (\ref{nor1}) becomes%
\begin{equation*}
g_{j}(x_{1},x_{2})=\lambda x_{j}+o(|x|^{M}),j=1,2,
\end{equation*}%
and then the $M$ th iteration $g^{M}=(g_{1}^{(M)},g_{2}^{(M)})$ has the form
of
\begin{equation*}
\begin{array}{c}
g_{1}^{(M)}(x_{1},x_{2})=x_{1}+o(|x|^{M}), \\
g_{2}^{(M)}(x_{1},x_{2})=x_{2}+o(|x|^{M}).%
\end{array}%
\end{equation*}%
Then, by Cronin's theorem and Lemma \ref{lem2-4}, we have%
\begin{equation*}
\mu _{f^{M}}(0)=\mu _{g^{M}}(0)=\pi _{g^{M}-id}(0)\geq (M+1)^{2},
\end{equation*}%
and then, by (\ref{all}) and by Lemmas \ref{lem2-5} and \ref{lem2-3} (ii),
we have
\begin{equation*}
P_{M}(f,0)=\mu _{f^{M}}(0)-P_{1}(f,0)\geq (M+1)^{2}-1>2M,
\end{equation*}%
and then, by Corollary \ref{do1}, $\mathcal{O}_{M}(f,0)>2.$ This completes
the proof.\medskip
\end{proof}

\begin{description}
\item[\textbf{Part 2}] \textbf{(b2) }$\Rightarrow $\textbf{\ (\ref{tar}%
).\medskip }
\end{description}

\begin{proof}
We first show that%
\begin{equation}
\mu _{f^{M}}(0)>M+1.  \label{kkk}
\end{equation}

By Lemma \ref{nor} and Corollary \ref{tra-d0}, there exists a polynomial
transform $(y_{1},y_{2})=H(x_{1},x_{2})$ in the form of (\ref{tra}), such
that the $M$ th iteration $g^{M}=(g_{1}^{(M)},g_{2}^{(M)})$ of the germ $%
g=(g_{1},g_{2})=H^{-1}\circ f\circ H$ has the expression%
\begin{equation*}
\begin{array}{c}
g_{1}^{(M)}(x_{1},x_{2})=\lambda
_{1}^{M}x_{1}+%
\sum_{_{i_{1}+i_{2}=2}}^{2M}C_{i_{1}i_{2}}^{1}x_{1}^{i_{1}}x_{2}^{i_{2}}+o(|x|^{2M}),
\\
g_{2}^{(M)}(x_{1},x_{2})=\lambda
_{2}^{M}x_{2}+%
\sum_{_{i_{1}+i_{2}=2}}^{2M}C_{i_{1}i_{2}}^{2}x_{1}^{i_{1}}x_{2}^{i_{2}}+o(|x|^{2M}),%
\end{array}%
\end{equation*}%
in a neighborhood of the origin, where, for $j=1$ and $2,$
\begin{equation}
C_{i_{1}i_{2}}^{j}\neq 0\mathrm{\ only\ if\ }\lambda _{j}=\lambda
_{1}^{i_{1}}\lambda _{2}^{i_{2}}.  \label{g1}
\end{equation}

By (b2), both $\lambda _{1}$ and $\lambda _{2}$ are primitive $M$ th roots
of unity. Thus, by (\ref{g1}) and (b2), $C_{i_{1}0}^{1}\neq 0$ only if $%
M|(i_{1}-1),$ and $C_{0i_{2}}^{1}\neq 0$ only if $i_{2}\geq \alpha ;$ $%
C_{0i_{2}}^{2}\neq 0$ only if $M|(i_{2}-1),$ and $C_{i_{1}0}^{2}\neq 0$ only
if $i_{1}\geq \beta .$ Thus, considering that $\lambda _{1}^{M}=\lambda
_{2}^{M}=1$ and that
\begin{equation*}
o(|x|^{2M})=x_{1}^{M+1}O(1)+x_{2}^{M+1}O(1),
\end{equation*}%
we can write%
\begin{equation}
\left\{
\begin{array}{c}
g_{1}^{(M)}(x_{1},x_{2})=x_{1}+x_{1}^{M+1}O(1)+O(1)x_{1}%
\sum_{i_{1}+i_{2}=1}^{2M}D_{i_{1}i_{2}}^{1}x_{1}^{i_{1}}x_{2}^{i_{2}}+x_{2}^{\alpha }O(1),
\\
g_{2}^{(M)}(x_{1},x_{2})=x_{2}+x_{1}^{\beta
}O(1)+O(1)x_{2}%
\sum_{i_{1}+i_{2}=1}^{2M}D_{i_{1}i_{2}}^{2}x_{1}^{i_{1}}x_{2}^{i_{2}}+x_{2}^{M+1}O(1),%
\end{array}%
\right.  \label{nor3}
\end{equation}%
in a neighborhood of the origin$,$ where
\begin{equation}
D_{i_{1}i_{2}}^{j}\neq 0\ \mathrm{only\ if\ }1=\lambda _{1}^{i_{1}}\lambda
_{2}^{i_{2}},j=1,2.  \label{g2}
\end{equation}

By condition (b2), we may write
\begin{equation*}
M=\beta m+\gamma ,
\end{equation*}%
where $m$ and $\gamma $ are positive integers with $\gamma <\beta .$ We
first assume%
\begin{equation}
\lambda _{1}=e^{\frac{2\pi i}{\beta m+\gamma }},\lambda _{2}=e^{\frac{2\beta
\pi i}{\beta m+\gamma }}.  \label{aid-2}
\end{equation}%
Then, for any pair $(i_{1},i_{2})$ of nonnegative integers with $\lambda
_{1}^{i_{1}}\lambda _{2}^{i_{2}}=1$ we have
\begin{equation}
i_{1}+\beta i_{2}=0(\mathrm{mod}M\mathrm{),}  \label{g3}
\end{equation}%
and then by (\ref{g2}), (\ref{g3}) and (\ref{aid}), putting $%
(x_{1},x_{2})=h(z_{1},z_{2})=(z_{1},z_{2}^{\beta }),$ we conclude that the
two components of the germ $(g^{M}-id)\circ h$ have the following
expression.
\begin{eqnarray*}
g_{1}^{(M)}\circ h(z_{1},z_{2})-z_{1}
&=&z_{1}^{M+1}O(1)+O(1)z_{1}%
\sum_{i_{1}+i_{2}=2}^{M}D_{i_{1}i_{2}}^{1}z_{1}^{i_{1}}z_{2}^{\beta
i_{2}}+z_{2}^{\alpha \beta }O(1) \\
&=&z_{1}\phi _{M}(z_{1},z_{2})+\mathrm{higher\ terms}, \\
g_{2}^{(M)}\circ h(z_{1},z_{2})-z_{2}^{\beta } &=&z_{1}^{\beta
}O(1)+O(1)z_{2}^{\beta
}\sum_{i_{1}+i_{2}=2}^{M}D_{i_{1}i_{2}}^{2}z_{1}^{i_{1}}z_{2}^{\beta
i_{2}}+z_{2}^{\beta (M+1)}O(1) \\
&=&az_{1}^{\beta }+\mathrm{higher\ terms},
\end{eqnarray*}%
where $\phi _{M}$ is a homogeneous polynomial of degree $M,$ in $z_{1}$ and $%
z_{2},$ and $a$ is a constant. In other words,%
\begin{equation}
\left\{
\begin{array}{l}
g_{1}^{(M)}\circ h(z_{1},z_{2})-z_{1}=z_{1}\phi _{M}(z_{1},z_{2})+\mathrm{%
higher\ terms}, \\
g_{2}^{(M)}\circ h(z_{1},z_{2})-z_{2}^{\beta }=az_{1}^{\beta }+\mathrm{%
higher\ terms}.%
\end{array}%
\right.   \label{ddd}
\end{equation}

Since we have assumed that the origin is an isolated fixed point of $f^{M},$
it is an isolated fixed point of $g^{M}$ by Lemma \ref{lem2-4}, and then the
origin is an isolated zero of the germ $(g^{M}-id)\circ h\in \mathcal{O}(%
\mathbb{C}^{2},0,0),$ which has the expression (\ref{ddd}). Thus, by
Cronin's theorem, the zero order $\pi _{(g^{M}-id)\circ h}(0)$ of the germ $%
(g^{M}-id)\circ h$ at the origin is not smaller than $\beta (M+1),$ and
since the origin is not an isolated solution of the system of equations
\begin{eqnarray*}
z_{1}\phi _{M}(z_{1},z_{2}) &=&0, \\
az_{1}^{\beta } &=&0,
\end{eqnarray*}%
we have by Cronin's theorem,%
\begin{equation*}
\pi _{(g^{M}-id)\circ h}(0)>\beta (M+1).
\end{equation*}%
Thus we have, by the fact $\pi _{h}(0)=\beta $ and Lemma \ref{prod-ind},
that
\begin{equation*}
\mu _{g^{M}}(0)=\pi _{g^{M}-id}(0)=\pi _{(g^{M}-id)\circ h}(0)/\pi
_{h}(0)>M+1.
\end{equation*}

Thus, by Lemma \ref{lem2-4}, we have proved (\ref{kkk}) under the assumption
(\ref{aid-2}).

When (\ref{aid-2}) fails, we first show that there exists a positive integer
$d$ such that $\Lambda _{1}=\lambda _{1}^{d}$ and $\Lambda _{2}=\lambda
_{2}^{d}$ has expression (\ref{aid-2}), say%
\begin{equation}
\Lambda _{1}=\lambda _{1}^{d}=e^{\frac{2\pi i}{\beta m+\gamma }},\Lambda
_{2}=\lambda _{2}^{d}=e^{\frac{2\beta \pi i}{\beta m+\gamma }}.
\label{aid-1}
\end{equation}

Since $\lambda _{1}$ is a primitive $M$ th root of unity, we may assume $%
\lambda _{1}=e^{\frac{2k_{1}\pi i}{M}},$ where $k_{1}$ is a positive integer
with $(k_{1},M)=1.$ Then there exists a positive integer $d$ such that
\begin{equation}
dk_{1}+cM=1,  \label{aid1}
\end{equation}%
where $c$ is an integer. Then,
\begin{equation}
\Lambda _{1}=\lambda _{1}^{d}=e^{\frac{2dk_{1}\pi i}{M}}=e^{\frac{2\pi i}{M}%
}.  \label{j5}
\end{equation}%
By (\ref{aid1}), one has $(d,M)=1,$ and then
\begin{equation*}
\Lambda _{2}=\lambda _{2}^{d}
\end{equation*}%
is still a primitive $M$ th root of unity, since $\lambda _{2}$ is a
primitive $M$ th root of unity.

Now, both $\Lambda _{1}$ and $\Lambda _{2}$ are primitive $M$ th roots of
unity, and by (\ref{aid})$,$ we still have%
\begin{eqnarray*}
\Lambda _{1}^{\beta } &=&\lambda _{1}^{d\beta }=\lambda _{2}^{d}=\Lambda
_{2}, \\
\Lambda _{2}^{\alpha } &=&\lambda _{2}^{d\alpha }=\lambda _{1}^{d}=\Lambda
_{1},
\end{eqnarray*}%
and then by (\ref{j5}), we have (\ref{aid-1}) (recall that $M=\beta m+\gamma
$).

Let $F=f^{d},$ the $d$ th iteration of $f$. Then $\Lambda _{1}$ and $\Lambda
_{2}$ are the two eigenvalues of $DF(0).$ Thus, the above argument for the
case (\ref{aid-2}) works for $\Lambda _{1}$, $\Lambda _{2}$ and $F,$
provided that $0$ is an isolated fixed point of $F^{M}=f^{dM}.$ On the other
hand, it is clear that the two eigenvalues of $Df^{M}(0)$ are both equal to $%
1,$ and therefore, by the assumption that $0$ is an isolated fixed point of $%
f^{M}$ and by Lemma \ref{SS-1}, $0$ is an isolated fixed point of $%
F^{M}=\left( f^{M}\right) ^{d}.$ Thus, applying the above argument to $%
\Lambda _{1},\Lambda _{2}$ and $F,$ we have%
\begin{equation*}
\mu _{F^{M}}(0)>M+1.
\end{equation*}%
and again by Lemma \ref{SS-1} we have%
\begin{equation*}
\mu _{f^{M}}(0)=\mu _{\left( f^{M}\right) ^{d}}(0)=\mu _{F^{M}}(0)>M+1,
\end{equation*}%
say, (\ref{kkk}) holds, and then we have proved (\ref{kkk}) completely.

Thus, by (\ref{all}), (\ref{kkk}) and by Lemmas \ref{lem2-5} and \ref{lem2-3}
(ii), we have
\begin{equation*}
P_{M}(f,0)=\mu _{f^{M}}(0)-P_{1}(f,0)>M,
\end{equation*}%
and then,
\begin{equation*}
\mathcal{O}_{M}(f,0)=P_{M}(f,0)/M\geq 2,
\end{equation*}%
for $\mathcal{O}_{M}(f,0)$ is an integer. This completes the proof.\medskip
\end{proof}

\begin{description}
\item[\textbf{Part 3}] \textbf{(b3) }$\Rightarrow $\textbf{\ (\ref{tar}).}%
\medskip
\end{description}

\begin{proof}
In case (b3), there exists an integer $d>1$ such that $M=m_{2}=dm_{1}$. Then
$\lambda _{1}$ is a primitive $m_{1}$ th root of unity and $\lambda _{2}$ is
a primitive $dm_{1}$ th root of unity, and%
\begin{equation}
\lambda _{2}^{d}\neq \lambda _{1}.  \label{ee}
\end{equation}%
We first prove the following two conclusions.

\textrm{(i)} If $j_{1}$ and $j_{2}$ are integers with $\lambda
_{1}^{j_{1}}\lambda _{2}^{j_{2}}=1$, then $d|j_{2}.$

\textrm{(ii)} For any positive integer $d_{1}$ with $\lambda
_{2}^{d_{1}}=\lambda _{1},$%
\begin{equation*}
d_{1}\geq 2d.
\end{equation*}

We may assume $\lambda _{1}=e^{\frac{2k_{1}\pi i}{m_{1}}}$ and $\lambda
_{2}=e^{\frac{2k_{2}\pi i}{dm_{1}}},$ where $k_{1}$ and $k_{2}$ are positive
integers with $k_{1}<m_{1}$, $k_{2}<dm_{1}$ and
\begin{equation}
(k_{1},m_{1})=(k_{2},dm_{1})=1.  \label{mon}
\end{equation}

If $j_{1}$ and $j_{2}$ are integers with $\lambda _{1}^{j_{1}}\lambda
_{2}^{j_{2}}=1$, then we have%
\begin{equation*}
\frac{j_{1}k_{1}}{m_{1}}+\frac{j_{2}k_{2}}{dm_{1}}=0(\mathrm{mod}1),
\end{equation*}%
and then%
\begin{equation*}
j_{1}k_{1}d+j_{2}k_{2}=0(\mathrm{mod}\left( dm_{1}\right) ).
\end{equation*}%
Thus by (\ref{mon}), we have $d|j_{2},$ and (i) is proved. If $\lambda
_{2}^{d_{1}}=\lambda _{1},$ then by (\ref{ee}), $d\neq d_{1},$ and by (i), $%
d|d_{1},$ which implies (ii).

Next, we show that there exists a polynomial transform $%
(y_{1},y_{2})=H(x_{1},x_{2})$ in the form of (\ref{tra}), such that for each
positive integer $k,$ the $k$ th iteration $g^{k}=(g_{1}^{(k)},g_{2}^{(k)})$
of the germ
\begin{equation}
g=(g_{1},g_{2})=H^{-1}\circ f\circ H  \label{dddd}
\end{equation}%
has the expression%
\begin{equation}
\left\{
\begin{array}{l}
g_{1}^{(k)}(x_{1},x_{2})=\lambda _{1}^{k}x_{1}+x_{1}^{m_{1}+1}(a^{(k)}+o(1))
\\
\ \ \ \ \ \ \ \ \ \ \ \ \ \ \ \ \ \
+x_{1}^{2}x_{2}^{d}O(1)+x_{2}^{2d}O(1),\medskip \\
g_{2}^{(k)}(x_{1},x_{2})=\lambda
_{2}^{k}x_{2}+x_{1}^{m_{1}}x_{2}O(1)+x_{1}x_{2}^{d+1}O(1) \\
\ \ \ \ \ \ \ \ \ \ \ \ \ \ \ \ \ \ +x_{2}^{2d+1}O(1)+x_{1}^{dm_{1}+1}o(1),%
\end{array}%
\right.  \label{nor4}
\end{equation}%
where $a^{(k)}$ is a constant for each $k.$

By Lemma \ref{nor} and Corollary \ref{tra-d0}, there exists a polynomial
transform $H$ in the form of (\ref{tra}) such that the $k$ th iteration $%
g^{k}=(g_{1}^{(k)},g_{2}^{(k)})$ of the germ (\ref{dddd}) has the expression%
\begin{equation}
g_{j}^{(k)}(x_{1},x_{2})=\lambda
_{j}^{k}x_{j}+%
\sum_{_{i_{1}+i_{2}=2}}^{3m_{1}d}C_{i_{1}i_{2}}^{kj}x_{1}^{i_{1}}x_{2}^{i_{2}}+o(|x|^{3dm_{1}})%
\mathrm{,\;}j=1,2,  \label{aa1}
\end{equation}%
in a neighborhood of the origin, in which $i_{1}$ and $i_{2}$ are
nonnegative integers and for $j=1$ and $2,$
\begin{equation}
C_{i_{1}i_{2}}^{kj}\neq 0\mathrm{\ only\ if\ }\lambda _{j}=\lambda
_{1}^{i_{1}}\lambda _{2}^{i_{2}}.  \label{aaa}
\end{equation}

By (\ref{aaa}), for $j=1$ and any pair $(i_{1},i_{2})$ in the sum with $%
C_{i_{1}i_{2}}^{k1}\neq 0,$ we have $\lambda _{1}^{i_{1}-1}\lambda
_{2}^{i_{2}}=1,$ and then by the assumption that $\lambda _{1}$ and $\lambda
_{2}$ are primitive $m_{1}$ th and $m_{2}$ th roots of unity (note that $%
m_{2}=dm_{1})$, respectively, for such pair $(i_{1},i_{2}),$ the following
conclusions from (1) to (4) hold (note that $i_{1}+i_{2}\geq 2$).

(1) If $i_{1}=0,$ then $i_{2}\geq 2d$ (by (\ref{ee}) and (ii))$;$

(2) If $i_{2}=0,$ then $i_{1}\geq m_{1}+1$;

(3) If $i_{1}=1,$ then $i_{2}\geq dm_{1}\geq 2d;$

(4) If $i_{1}\geq 2$ and $i_{2}\geq 1$, then $i_{2}\geq d$ (by (i))$.$

On the other hand, since $d>1$ and $m_{1}>1,$ it is clear that any term of $%
o(|x|^{3dm_{1}})$ has either a factor $x_{1}^{dm_{1}+2},$ or a factor $%
x_{2}^{2d+1},$ and then one can write%
\begin{equation}
o(|x|^{3dm_{1}})=x_{1}^{dm_{1}+2}O(1)+x_{2}^{2d+1}O(1).  \label{aa2}
\end{equation}%
Thus, by (1)--(4) we obtain the first equation in (\ref{nor4}).

By (\ref{aaa}), for $j=2$ and any pair $(i_{1},i_{2})$ with $%
C_{i_{1}i_{2}}^{k2}\neq 0,$ we have $\lambda _{1}^{i_{1}}\lambda
_{2}^{i_{2}-1}=1,$ and then by the assumption, for such pair $(i_{1},i_{2})$
the following conclusions from (5) to (8) hold (note that $i_{1}+i_{2}\geq 2$%
).

(5) $i_{2}\neq 0;$

(6) If $i_{1}=0,$ then $i_{2}\geq dm_{1}+1\geq 2d+1$;

(7) If $i_{2}=1,$ then $i_{1}\geq m_{1};$

(8) If $m_{1}>i_{1}\geq 1$ and $i_{2}\geq 2,\;$then $i_{2}\geq d+1$ (by (i))$%
.$

Thus, by (\ref{aa2}), the second equation in (\ref{nor4}) holds, and (\ref%
{nor4}) is proved.

Since the origin is an isolated fixed point of $f^{dm_{1}},$ it is an
isolated fixed point of $f^{m_{1}},g^{m_{1}}$ and $g^{dm_{1}}$ as well, by
Lemma \ref{lem2-4}. To complete the proof, we first show (\ref{tar}) under
the assumption that%
\begin{equation*}
a^{(1)}\neq 0.
\end{equation*}%
Then it is easy to see that the coefficient $a^{(k)}$ in (\ref{nor4})
satisfies
\begin{equation}
a^{(k)}=k\lambda _{1}^{k-1}a^{(1)}\neq 0,k\in \mathbb{N}.  \label{nor5}
\end{equation}

For $k=m_{1},$ by (\ref{nor4}) and by the fact $\lambda _{1}^{m_{1}}=1\ $and
$\lambda _{2}^{m_{1}}\neq 1,$ one can write%
\begin{eqnarray*}
g_{1}^{(m_{1})}(x_{1},x_{2})-x_{1}
&=&x_{1}^{m_{1}+1}(a^{(m_{1})}+o(1))+x_{2}o(1), \\
g_{2}^{(m_{1})}(x_{1},x_{2})-x_{2} &=&bx_{2}+x_{2}o(1)+x_{1}^{m_{1}+1}o(1).
\end{eqnarray*}%
where $b=\lambda _{2}^{m_{1}}-1\neq 0,$ and then by (\ref{nor5}) and
Corollary \ref{c2} we have%
\begin{equation}
\mu _{g^{m_{1}}}(0)=\pi _{g^{m_{1}}-id}(0)=m_{1}+1.  \label{nor5+1}
\end{equation}

For\label{?2} $k=dm_{1},$ by (\ref{nor4}) we have
\begin{equation}
\left\{
\begin{array}{l}
g_{1}^{(dm_{1})}(x_{1},x_{2})=x_{1}+x_{1}^{m_{1}+1}O(1)+x_{1}^{2}x_{2}^{d}O(1)+x_{2}^{2d}O(1),\medskip
\\
g_{2}^{(dm_{1})}(x_{1},x_{2})=x_{2}+x_{1}^{m_{1}}x_{2}O(1)+x_{1}x_{2}^{d+1}O(1)
\\
\ \ \ \ \ \ \ \ \ \ \ \ \ \ \ \ \ \ \ \ \ \ \
+x_{2}^{2d+1}O(1)+x_{1}^{dm_{1}+1}o(1),%
\end{array}%
\right.  \label{c62}
\end{equation}%
and then, putting $(x_{1},x_{2})=h(z_{1},z_{2})=(z_{1}^{d},z_{2}^{m_{1}}),$
we conclude that the germ $G=(G_{1},G_{2})=\left( g^{dm_{1}}-id\right) \circ
h$ has the expression%
\begin{eqnarray*}
G_{1}(z_{1},z_{2}) &=&\phi (z_{1},z_{2})+\mathrm{higher\ terms,} \\
G_{2}(z_{1},z_{2}) &=&\psi (z_{1},z_{2})+\mathrm{higher\ terms,}
\end{eqnarray*}%
where $\phi $ and $\psi $ are homogeneous polynomials of degrees\label{?3} $%
dm_{1}+d$ and $dm_{1}+m_{1},$ respectively (note that $d>1$ and $m_{1}>1)$.
Since the origin is an isolated fixed point of $f^{M},$ it is also an
isolated fixed point of $g^{M}$ by Lemma \ref{lem2-4}, and then it is an
isolated zero of $G.$ Therefore, by Cronin's theorem we have%
\begin{equation*}
\pi _{G}(0)\geq (dm_{1}+d)(dm_{1}+m_{1})=dm_{1}(m_{1}+1)(d+1),
\end{equation*}%
and then, we have $\mu _{g^{dm_{1}}}(0)=\pi _{G}(0)/\pi _{h}(0)\geq
(m_{1}+1)(d+1),$ say,%
\begin{equation}
\mu _{g^{dm_{1}}}(0)>dm_{1}+m_{1}+1.  \label{c64}
\end{equation}

Thus, we have by (\ref{nor5+1}), (\ref{c64}) and Lemma \ref{lem2-4} that%
\begin{equation}
\mu _{f^{m_{1}}}(0)=m_{1}+1\ \mathrm{and\ }\mu
_{f^{dm_{1}}}(0)>dm_{1}+m_{1}+1.  \label{j11}
\end{equation}%
By Lemma \ref{lem2-5}, each periodic point of the linear part of $f$ at the
origin has period $1,$ $m_{1}$ or $dm_{1}.$ Thus, by Lemma \ref{lem2-3}
(ii), we have%
\begin{eqnarray*}
P_{dm_{1}}(f,0) &=&\mu _{f^{dm_{1}}}(0)-P_{m_{1}}(f,0)-P_{1}(f,0), \\
P_{m_{1}}(f,0) &=&\mu _{f^{m_{1}}}(0)-P_{1}(f,0),
\end{eqnarray*}%
and then, we have by (\ref{j11}) that%
\begin{equation*}
P_{dm_{1}}(f,0)=\mu _{f^{dm_{1}}}(0)-\mu _{f^{m_{1}}}(0)>dm_{1},
\end{equation*}%
and then $\mathcal{O}_{dm_{1}}(f,0)=P_{dm_{1}}(f,0)/(dm_{1})>1.$ But $%
\mathcal{O}_{dm_{1}}(f,0)$ is an integer, we have $\mathcal{O}%
_{dm_{1}}(f,0)\geq 2,$ and we have proved (\ref{tar}) under the assumption $%
a^{(1)}\neq 0.$

If $a^{(1)}=0,$ then for $\varepsilon \neq 0$ consider the mapping $%
g_{\varepsilon }=(g_{1,\varepsilon },g_{2,\varepsilon })$ given by%
\begin{equation}
\left\{
\begin{array}{l}
g_{1,\varepsilon }(x_{1},x_{2})=g_{1}(x_{1},x_{2})+\varepsilon
x_{1}^{m_{1}+1}, \\
g_{2,\varepsilon }(x_{1},x_{2})=g_{2}(x_{1},x_{2}),%
\end{array}%
\right.  \notag
\end{equation}%
which is obtain from (\ref{nor4}) with $k=1$, by just replacing $a^{(1)},$
the coefficient of $x_{1}^{m_{1}+1}$ of the power series of $g_{1}^{(1)},$
with $\varepsilon $ (note that $%
g=g^{1}=(g_{1},g_{2})=(g_{1}^{(1)},g_{2}^{(1)})$). Then $f_{\varepsilon
}=H\circ g_{\varepsilon }\circ H^{-1}$ converges to $f$ uniformly in a
neighborhood of the origin as $\varepsilon \rightarrow 0,$ and then, $0$ is
also an isolated fixed point of $f_{\varepsilon }^{dm_{1}}$ for sufficiently
small $\varepsilon ,$ by the assumption that $0$ is an isolated fixed point
of $f^{dm_{1}}$ and Lemma \ref{lem2-1}.

On the other hand, $g_{\varepsilon }=H^{-1}\circ f_{\varepsilon }\circ H$ is
in the form of (\ref{aa1}) for $k=1,$ together (\ref{aaa}). Thus, by
Corollary \ref{tra-d0}, the $k$ th iteration $g_{\varepsilon }^{k}$ is still
in the form of (\ref{aa1})$,$ together (\ref{aaa}). Thus, repeating the
process for proving (\ref{nor4}), we can prove that $g_{\varepsilon }^{k}$
is still in the form of (\ref{nor4}), more precisely%
\begin{eqnarray*}
g_{1}^{(k)}(x_{1},x_{2}) &=&\lambda
_{1}^{k}x_{1}+x_{1}^{m_{1}+1}(A^{(k)}+o(1)) \\
&&+x_{1}^{2}x_{2}^{d}O(1)+x_{2}^{2d}O(1), \\
g_{2}^{(k)}(x_{1},x_{2}) &=&\lambda
_{2}^{k}x_{2}+x_{1}^{m_{1}}x_{2}O(1)+x_{1}x_{2}^{d+1}O(1) \\
&&+x_{2}^{2d+1}O(1)+x_{1}^{dm_{1}+1}o(1).
\end{eqnarray*}%
But here $A^{(1)}=\varepsilon \neq 0.$

Therefore, all the above arguments for the case $a^{(1)}\neq 0$ apply to $%
f_{\varepsilon },g_{\varepsilon }=H^{-1}\circ f_{\varepsilon }\circ H$ and $%
g_{\varepsilon }^{k},$ and then we have $\mathcal{O}_{dm_{1}}(f_{\varepsilon
},0)\geq 2.$ Hence, by Lemma \ref{lem2-2}, we have (\ref{tar}) in the case $%
a^{(1)}=0$ as well. This completes the proof.\medskip
\end{proof}

\begin{description}
\item[\textbf{Part 4}] \textbf{(b4) }$\Rightarrow $\textbf{\ (\ref{tar}).}%
\medskip
\end{description}

\begin{proof}
By (b4), there exist positive integers $d>1,$ $n_{1}>1$ and $n_{2}>1$ such
that $n_{1}$ and $n_{2}$ are relatively prime and
\begin{equation*}
m_{1}=dn_{1},\ m_{2}=dn_{2},\ M=dn_{1}n_{2}.
\end{equation*}%
Then, the two eigenvalues $\lambda _{1}$ and $\lambda _{2}$ of $Df(0)$ are
primitive $dn_{1}$-th, and $dn_{2}$-th roots of unity, respectively.

We first show that there exists a polynomial transform $H$ in the form of (%
\ref{tra}) in a neighborhood of the origin such that for each $k\in \mathbb{N%
},$ the $k$ th iteration $g^{k}=(g_{1}^{(k)},g_{2}^{(k)})$ of the germ
\begin{equation}
g=H^{-1}\circ f\circ H=(g_{1},g_{2})  \label{h10}
\end{equation}%
has the expression%
\begin{equation}
\left\{
\begin{array}{lll}
g_{1}^{(k)}(x_{1},x_{2}) & = & \lambda
_{1}^{k}x_{1}+x_{1}^{dn_{1}+1}(a^{(k)}+o(1))+x_{1}^{n_{1}+1}x_{2}^{n_{2}}O(1)
\\
&  & +x_{1}x_{2}^{dn_{2}}O(1)+x_{2}^{2dn_{2}+1}O(1),\medskip \\
g_{2}^{(k)}(x_{1},x_{2}) & = & \lambda
_{2}^{k}x_{2}+x_{1}^{dn_{1}}x_{2}O(1)+x_{1}^{n_{1}}x_{2}^{n_{2}+1}O(1) \\
&  & +x_{2}^{dn_{2}+1}(b^{(k)}+o(1))+x_{1}^{2dn_{1}+1}O(1),%
\end{array}%
\right.  \label{tra2}
\end{equation}%
where $a^{(k)}$ and $b^{(k)}$ are constants.

By Lemma \ref{nor} and Corollary \ref{tra-d0}, there exists a polynomial
coordinate transform $H$ in the form of (\ref{tra}) in a neighborhood of the
origin, such that each component of the $k$ th iteration $%
g^{k}=(g_{1}^{(k)},g_{2}^{(k)})$ of the germ (\ref{h10}) has the expression%
\begin{equation}
g_{j}^{(k)}(x_{1},x_{2})=\lambda
_{j}^{k}x_{j}+%
\sum_{_{i_{1}+i_{2}=2}}^{4dn_{1}n_{2}}C_{i_{1}i_{2}}^{kj}x_{1}^{i_{1}}x_{2}^{i_{2}}+o(|x|^{4dn_{1}n_{2}}),j=1,2,
\label{ad-4.1.2}
\end{equation}%
in a neighborhood of the origin, where for each $j=1,2,$ the sum in (\ref%
{ad-4.1.2}) extends over all $2$-tuples $(i_{1},i_{2})$ of nonnegative
integers with%
\begin{equation}
2\leq i_{1}+i_{2}\leq 4dn_{1}n_{2}  \label{0}
\end{equation}%
and%
\begin{equation}
\lambda _{j}=\lambda _{1}^{i_{1}}\lambda _{2}^{i_{2}}.  \label{1}
\end{equation}

We may write
\begin{equation*}
\lambda _{j}=e^{\frac{2k_{j}\pi i}{dn_{j}}},j=1,2,
\end{equation*}%
where $k_{j}$ is a positive integer such that $(k_{j},dn_{j})=1,j=1,2.$
Recall that $(k_{j},dn_{j})$ denotes the largest common divisor of $k_{j}$
and $dn_{j}.$

Now, first assume $j=1$ and let $(i_{1},i_{2})$ be any $2$-tuple that
satisfies (\ref{0}) and (\ref{1}). Then we have%
\begin{equation*}
\frac{(i_{1}-1)k_{1}}{dn_{1}}+\frac{i_{2}k_{2}}{dn_{2}}=0(\mathrm{{mod}1),}
\end{equation*}%
and then%
\begin{equation}
(i_{1}-1)k_{1}n_{2}+i_{2}k_{2}n_{1}=0(\mathrm{mod}\left( dn_{1}n_{2}\right)
).  \label{3}
\end{equation}%
This implies that $n_{1}|(i_{1}-1)$ and $n_{2}|i_{2},$ since
\begin{equation*}
(k_{1},n_{1})=(k_{2},n_{2})=(n_{1},n_{2})=1.
\end{equation*}%
Therefore, for $j=1,$ each nonzero term $%
C_{i_{1}i_{2}}^{k1}x_{1}^{i_{1}}x_{2}^{i_{2}}$ in (\ref{ad-4.1.2}) is of
type $x_{1}^{sn_{1}+1}x_{2}^{tn_{2}}O(1),$ in which $s=(i_{1}-1)/n_{1}$ and $%
t=i_{2}/n_{2}$ are integers. In particular, when $i_{2}=0,$ by (\ref{3}) we
have $\left( dn_{1}\right) |\left( i_{1}-1\right) ,$ and then $i_{1}\geq
dn_{1}+1$ by (\ref{0}); and when $i_{1}=1,$ by (\ref{3}) we have $\left(
dn_{2}\right) |i_{2},$ and then $i_{2}\geq dn_{2}$ by (\ref{0}). On the
other hand, each term in the power series $o(|x|^{4dn_{1}n_{2}})$ has the
form $cx_{1}^{i_{1}}x_{2}^{i_{2}}$ with $i_{1}+i_{2}>4dn_{1}n_{2}.$ Thus
each term in $o(|x|^{4dn_{1}n_{2}})$ has either a factor $x_{1}^{dn_{1}+2}$
or a factor $x_{2}^{2dn_{2}+1}.$ Therefore, the first equality of (\ref{tra2}%
) holds for some constant $a^{(k)}$.

For the same reason, the second equality in (\ref{tra2}) holds for some
constant $b^{(k)}$.

By the assumption, we have $\lambda _{1}^{dn_{1}n_{2}}=\lambda
_{2}^{dn_{1}n_{2}}=1,$ and then, by taking $k=dn_{1}n_{2}$ in (\ref{tra2}),
we have%
\begin{eqnarray*}
g_{1}^{(dn_{1}n_{2})}(x_{1},x_{2})-x_{1}
&=&x_{1}^{dn_{1}+1}O(1)+x_{1}^{n_{1}+1}x_{2}^{n_{2}}O(1) \\
&&+x_{1}x_{2}^{dn_{2}}O(1)+x_{2}^{2dn_{2}+1}O(1), \\
g_{2}^{(dn_{1}n_{2})}(x_{1},x_{2})-x_{2}
&=&x_{1}^{dn_{1}}x_{2}O(1)+x_{1}^{n_{1}}x_{2}^{n_{2}+1}O(1) \\
&&+x_{2}^{dn_{2}+1}O(1)+x_{1}^{2dn_{1}+1}O(1).
\end{eqnarray*}%
By Lemma \ref{dd4}, the zero order $\pi _{g^{dn_{1}n_{2}}-id}(0)$ of $%
g^{dn_{1}n_{2}}-id$ at the origin is at least $1+dn_{1}+dn_{2}+2dn_{1}n_{2},$
in other words%
\begin{equation*}
\mu _{g^{dn_{1}n_{2}}}(0)\geq 1+dn_{1}+dn_{2}+2dn_{1}n_{2},
\end{equation*}%
which implies by Lemma \ref{lem2-4} that
\begin{equation}
\mu _{f^{dn_{1}n_{2}}}(0)\geq 1+dn_{1}+dn_{2}+2dn_{1}n_{2}.  \label{dd5}
\end{equation}

Now, let us first prove (\ref{tar}) under the assumption that, in the
expression (\ref{tra2}), $a^{(1)}\neq 0$ and $b^{(1)}\neq 0.$

Then it is easy to see that%
\begin{equation*}
a^{(dn_{1})}=dn_{1}\lambda _{1}^{dn_{1}-1}a^{(1)}\neq 0.
\end{equation*}%
Thus, by (\ref{tra2}) and the fact that $\lambda _{1}^{dn_{1}}=1$ and $%
\lambda _{2}^{dn_{1}}\neq 1$, we can write%
\begin{eqnarray*}
g_{1}^{(dn_{1})}(x_{1},x_{2})-x_{1}
&=&x_{1}^{dn_{1}+1}(a^{(dn_{1})}+o(1))+x_{2}o(1), \\
g_{2}^{(dn_{1})}(x_{1},x_{2})-x_{2} &=&cx_{2}+x_{2}o(1)+x_{1}^{dn_{1}+1}o(1),
\end{eqnarray*}%
where $c=\lambda _{2}^{dn_{1}}-1\neq 0,$ and then by Corollary \ref{c2} we
have%
\begin{equation*}
\mu _{g^{dn_{1}}}(0)=\pi _{g^{dn_{1}}-id}(0)=dn_{1}+1,
\end{equation*}%
and, repeating the above argument, by the assumption $b^{(1)}\neq 0,$ we have%
\begin{equation*}
\mu _{g^{dn_{2}}}(0)=dn_{2}+1.
\end{equation*}%
Thus by Lemma \ref{lem2-4} we have
\begin{equation}
\mu _{f^{dn_{1}}}(0)=dn_{1}+1,\mu _{f^{dn_{2}}}(0)=dn_{2}+1.  \label{ai+2}
\end{equation}%
On the other hand, by Lemmas \ref{lem2-3} (ii) and \ref{lem2-5}, we have%
\begin{eqnarray*}
\mu _{f^{dn_{1}}}(0) &=&P_{dn_{1}}(f,0)+P_{1}(f,0), \\
\mu _{f^{dn_{2}}}(0) &=&P_{dn_{2}}(f,0)+P_{1}(f,0), \\
\mu _{f^{dn_{1}n_{2}}}(0)
&=&P_{dn_{1}n_{2}}(f,0)+P_{dn_{1}}(f,0)+P_{dn_{2}}(f,0)+P_{1}(f,0).
\end{eqnarray*}%
Thus, by (\ref{all}) we have%
\begin{equation*}
P_{dn_{1}n_{2}}(f,0)=\mu _{f^{dn_{1}n_{2}}}-\mu _{f^{dn_{1}}}(0)-\mu
_{f^{dn_{2}}}(0)+1.
\end{equation*}%
and then, by (\ref{dd5}) and by (\ref{ai+2}), we have $P_{dn_{1}n_{2}}(f,0)%
\geq 2dn_{1}n_{2},$ and then (\ref{tar}) holds.

Now, we have proved (\ref{tar}) under the condition $a^{(1)}\neq 0$ and $%
b^{(1)}\neq 0.$ In general$,$ we consider $g_{\varepsilon
}=(g_{1,\varepsilon },g_{2,\varepsilon })\in \mathcal{O}(\mathbb{C}^{2},0,0)$
given by%
\begin{equation}
\left.
\begin{array}{c}
g_{1,\varepsilon }(x_{1},x_{2})=g_{1}(x_{1},x_{2})+\varepsilon x^{dn_{1}+1},
\\
g_{2,\varepsilon }(x_{1},x_{2})=g_{2}(x_{1},x_{2})+\varepsilon
x_{2}^{dn_{2}+1},%
\end{array}%
\right.  \label{2}
\end{equation}%
which is obtained from (\ref{tra2}) with $k=1$, by just replacing the
constants $a^{(1)}$ and $b^{(1)}$ with $\varepsilon $, and consider
\begin{equation*}
f_{\varepsilon }=H\circ g_{\varepsilon }\circ H^{-1}
\end{equation*}%
where $H$ is the transform in (\ref{h10}).

Since $f_{\varepsilon }$ uniformly converges to $f$ as $\varepsilon
\rightarrow 0$, for sufficiently small $\varepsilon ,$ the origin is an
isolated fixed point of $f_{\varepsilon }^{dn_{1}n_{2}}$ by Rouch\'{e}'s
theorem and the assumption that $0$ is an isolated fixed point of $%
f^{dn_{1}n_{2}}$ (note that $M=dn_{1}n_{2})$. Then, it is clear that for
sufficiently small $\varepsilon ,f_{\varepsilon }$ again satisfies condition
(b4) and $g_{\varepsilon }=H^{-1}\circ f_{\varepsilon }\circ H$ is in the
form of (\ref{ad-4.1.2}) of $k=1,$ together with (\ref{0}) and (\ref{1}),
thus $g_{\varepsilon }^{k}=(g_{1,\varepsilon }^{(k)},g_{2,\varepsilon
}^{(k)})$ is in the form of (\ref{ad-4.1.2}) for all $k\in \mathbb{N}$ by
Corollary \ref{tra-d0}, more precisely,
\begin{equation*}
\begin{array}{lll}
g_{1,\varepsilon }^{(k)}(x_{1},x_{2}) & = & \lambda
_{1}^{k}x_{1}+x_{1}^{dn_{1}+1}(A^{(k)}+o(1))+x_{1}^{n_{1}+1}x_{2}^{n_{2}}O(1)
\\
&  & +x_{1}x_{2}^{dn_{2}}O(1)+x_{2}^{2dn_{2}+1}O(1), \\
g_{2,\varepsilon }^{(k)}(x_{1},x_{2}) & = & \lambda
_{2}^{k}x_{2}+x_{1}^{dn_{1}}x_{2}O(1)+x_{1}^{n_{1}}x_{2}^{n_{2}+1}O(1) \\
&  & +x_{2}^{dn_{2}+1}(B^{(k)}+o(1))+x_{1}^{2dn_{1}+1}O(1),%
\end{array}%
\end{equation*}%
where $A^{(1)}=B^{(1)}=\varepsilon \neq 0$.

Thus, all the above arguments for the case $a^{(1)}b^{(1)}\neq 0$ apply to $%
f_{\varepsilon }$ and $g_{\varepsilon }^{k}=(g_{1,\varepsilon
}^{(k)},g_{2,\varepsilon }^{(k)}),$ and then we have $\mathcal{O}%
_{dn_{1}n_{2}}(f_{\varepsilon },0)\geq 2,$ and then by Lemma \ref{lem2-2},
we have $\mathcal{O}_{dn_{1}n_{2}}(f,0)\geq 2,$ say (\ref{tar}) holds again$%
. $ This completes the proof.
\end{proof}

\section{Proof of the Main Theorem: (A) $\Rightarrow $ (B)\label{PM}}

Assume that $M>1$ is a positive integer and assume that $A$ is a $2\times 2$
matrix such that the following condition holds.

\textrm{(C)} \emph{For any }$f\in \mathcal{O}(\mathbb{C}^{2},0,0)$\emph{\
such that }$Df(0)=A$\emph{\ and that }$0$\emph{\ is an isolated fixed point
of }$f^{M},$
\begin{equation}
\mathcal{O}_{M}(f,0)\geq 2.  \label{p1}
\end{equation}

We show that $A$ satisfies condition \textrm{(B)} in Theorem \ref{Th1}.

We first show that the two eigenvalues $\lambda _{1}$ and $\lambda _{2}$ of $%
A$ satisfy
\begin{equation}
\lambda _{1}^{M}=\lambda _{2}^{M}=1.  \label{p1+1}
\end{equation}

Assume (\ref{p1+1}) fails. Then we may assume $\lambda _{2}^{M}\neq 1,$ and
then by (C), Theorem \ref{Th0.1} and Lemma \ref{lem2-5}, the other
eigenvalue $\lambda _{1}$ must be a primitive $M$ th root of unity. But then
we shall obtain a contradiction by constructing a germ $f\in \mathcal{O}(%
\mathbb{C}^{2},0,0)$ such that $Df(0)=A$ and\ that $0$\ is an isolated fixed
point of $f^{M},$ but (\ref{p1}) fails. That is to say, (C) fails!

Since $\lambda _{1}$ is a primitive $M$ th root of unity and $\lambda
_{2}^{M}\neq 1,$ we have $\lambda _{1}\neq \lambda _{2},$ and then the
matrix $A$ is diagonalizable. Thus, by Lemma \ref{lem2-4}, to construct that
$f,$ we may assume that
\begin{equation*}
A=\left(
\begin{array}{cc}
\lambda _{1} & 0 \\
0 & \lambda _{2}%
\end{array}%
\right) .
\end{equation*}

We show that the germ $f(x_{1},x_{2})\in \mathcal{O}(\mathbb{C}^{2},0,0)$
given by%
\begin{equation*}
f(x_{1},x_{2})=(\lambda _{1}x_{1}+x_{1}^{M+1},\lambda _{2}x_{2})
\end{equation*}%
is the desired germ.

Since $\lambda _{1}$ is a primitive $M$ th root of unity and $M>1$, it is
clear that $\lambda _{1}\neq 1$ and that the $M$ th iteration $f^{M}$ of $f$
has the expression%
\begin{equation*}
f^{M}(x_{1},x_{2})=(x_{1}+cx_{1}^{M+1}(1+o(1)),\lambda _{2}^{M}x_{2}),
\end{equation*}%
where $c=M\lambda _{1}^{M-1}\neq 0.$ Then, by the assumption that $\lambda
_{2}^{M}\neq 1$ and by Cronin's theorem, the origin is an isolated zero of $%
f^{M}-id$ and the zero order is%
\begin{equation*}
\pi _{f^{M}-id}(0)=M+1,
\end{equation*}%
and then the origin is an isolated fixed point of $f^{M}$ with%
\begin{equation*}
\mu _{f^{M}}(0)=M+1.
\end{equation*}%
On the other hand, it is clear that $0$ is a simple fixed point of $f,$ and
then
\begin{equation*}
\mu _{f}(0)=P_{1}(f,0)=1.
\end{equation*}
Therefore, by Lemmas \ref{lem2-3} (ii) and \ref{lem2-5},
\begin{equation*}
P_{M}(f,0)=\mu _{f^{M}}(0)-P_{1}(f,0)=M,
\end{equation*}%
and then $\mathcal{O}_{M}(f,0)=1,$ say, (\ref{p1}) fails$.$ Hence, (\ref%
{p1+1}) holds.

By \textrm{(C)}, (\ref{p1+1}), Theorem \ref{Th0.1} and Lemma \ref{lem2-5},
we can conclude that there exists positive integers $m_{1}$ and $m_{2}$ such
that

(D) \emph{The two eigenvalues }$\lambda _{1}$\emph{\ and }$\lambda _{2}$%
\emph{\ of }$A$ \emph{are primitive }$m_{1}$\emph{\ th and }$m_{2}$\emph{\
th roots of unity, respectively, and}
\begin{equation}
M=[m_{1},m_{2}].  \label{p1+2}
\end{equation}

Without loss of generality, we assume that
\begin{equation*}
m_{1}\leq m_{2}\leq M.
\end{equation*}

If $m_{1}$ and $m_{2}$ do not satisfy any one of the conditions from (b1) to
(b4), then by (D) and Remark \ref{r1}, one of the following conditions must
be satisfied:

\textrm{(b1)'\ }$m_{1}=m_{2}=M,$ $\lambda _{1}=\lambda _{2}$ and $A$ is not
diagonalizable$.$

\textrm{(b2)'\ }$m_{1}=m_{2}=M$ and there exists positive integers $\alpha $
and $\beta \ $with $1<\alpha <M$ and $1<\beta <M$ such that%
\begin{equation}
\lambda _{1}^{\alpha }=\lambda _{2},\lambda _{2}^{\beta }=\lambda
_{1},\alpha \beta =M+1.  \label{p2}
\end{equation}

\textrm{(b3)'\ }$m_{2}=M,$ $m_{1}|m_{2},m_{1}<m_{2},$ and $\lambda
_{2}^{m_{2}/m_{1}}=\lambda _{1}.$

\textrm{(b4)'} $m_{1}$ and $m_{2}$ are relatively prime, and $m_{2}>m_{1}>1$.

We show that each condition from (b1)' to (b4)' contradicts condition (C).
This will be done in each case, by constructing a germ $F\in \mathcal{O}(%
\mathbb{C}^{2},0,0)$ such that $DF(0)=A$ and $0$ is an isolated fixed point
of $F^{M},$ but $\mathcal{O}_{M}(F,0)=1,$ say, (\ref{p1}) fails. We divide
this process into four parts.

\begin{description}
\item[\textbf{Part 1}] \textbf{(b1)' implies the existences of }$F$\textbf{.}
\end{description}

\begin{proof}
In case \textrm{(b1)'}, by Lemma \ref{lem2-4}, we may assume that%
\begin{equation*}
A=\left(
\begin{array}{cc}
\lambda _{1} & 0 \\
1 & \lambda _{1}%
\end{array}%
\right) .
\end{equation*}%
Then, consider the germ $F\in \mathcal{O}(\mathbb{C}^{2},0,0)$ given by%
\begin{equation*}
F(x_{1},x_{2})=(\lambda _{1}x_{1}+x_{2}^{M+1},x_{1}+\lambda _{1}x_{2}).
\end{equation*}%
Since $\lambda _{1}$ is a primitive $M$ th root of unity and $M>1,$ the
origin is a simple fixed point of $f$, say,%
\begin{equation}
P_{1}(F,0)=\mu _{F}(0)=1.  \label{pf1}
\end{equation}

It is easy to see by induction that the $M$ th iteration $F^{M}$ of $F$ has
the expression%
\begin{equation*}
\left( F^{M}(x_{1},x_{2})\right) ^{T}=\left(
\begin{array}{c}
x_{1}+x_{1}\phi _{M}(x_{1},x_{2})+M\lambda _{1}^{M-1}x_{2}^{M+1}+o(|x|^{M+1})
\\
M\lambda _{1}^{M-1}x_{1}+x_{2}+o(|x|)%
\end{array}%
\right) ,
\end{equation*}%
where $x=(x_{1},x_{2})$ and $\phi _{M}$ is a homogenous polynomial of degree
$M.$ Then $F^{M}(x_{1},x_{2})-(x_{1},x_{2})$ has the expression
\begin{equation*}
\left( F^{M}(x_{1},x_{2})-(x_{1},x_{2})\right) ^{T}=\left(
\begin{array}{c}
x_{1}\phi _{M}(x_{1},x_{2})+M\lambda _{1}^{M-1}x_{2}^{M+1}+o(|x|^{M+1}) \\
M\lambda _{1}^{M-1}x_{1}+o(|x|)%
\end{array}%
\right)
\end{equation*}%
and the origin is an isolated zero of the system%
\begin{eqnarray*}
x_{1}\phi _{M}(x_{1},x_{2})+M\lambda _{1}^{M-1}x_{2}^{M+1} &=&0, \\
M\lambda _{1}^{M-1}x_{1} &=&0.
\end{eqnarray*}%
Thus, by Cronin's theorem, the zero order of the germ $F^{M}-id$ at the
origin equals $M+1,$ and then
\begin{equation*}
\mu _{F^{M}}(0)=M+1.
\end{equation*}%
Thus, by (\ref{pf1}) and by Lemmas \ref{lem2-3} (ii) and \ref{lem2-5},
\begin{equation*}
P_{M}(F,0)=\mu _{F^{M}}(0)-P_{1}(F,0)=M,
\end{equation*}%
which implies $\mathcal{O}_{M}(F,0)=1.$
\end{proof}

In any other case, $A$ is diagonalizable, and by Lemma \ref{lem2-4}, we may
assume
\begin{equation*}
A=\left(
\begin{array}{cc}
\lambda _{1} & 0 \\
0 & \lambda _{2}%
\end{array}%
\right) .
\end{equation*}

\begin{description}
\item[\textbf{Part 2}] \textbf{(b2)'}$\ $\textbf{implies the existence of }$%
F $\textbf{.}
\end{description}

\begin{proof}
\textrm{In case (b2)',} consider the germ $F=(F_{1},F_{2})\in \mathcal{O}(%
\mathbb{C}^{2},0,0)$ given by%
\begin{eqnarray*}
F_{1}(x_{1},x_{2}) &=&\lambda _{1}x_{1}+x_{2}^{\beta }, \\
F_{2}(x_{1},x_{2}) &=&\lambda _{2}x_{2}+x_{1}^{\alpha }.
\end{eqnarray*}%
By (\textrm{D}) and (b2)', both $\lambda _{1}$ and $\lambda _{2}$ are
primitive $M$ th roots of unity, and then by (\ref{p2}), it is easy to see
that the $M$ th iteration $F^{M}=(F_{1}^{(M)},F_{2}^{(M)})$ of $F$ has the
expression%
\begin{eqnarray*}
F_{1}^{(M)}(x_{1},x_{2}) &=&x_{1}+M\lambda _{1}^{M-1}x_{2}^{\beta }+\mathrm{%
higher\ terms,} \\
F_{2}^{(M)}(x_{1},x_{2}) &=&x_{2}+M\lambda _{2}^{M-1}x_{1}^{\alpha }+\mathrm{%
higher\ terms;}
\end{eqnarray*}%
and then by Cronin's theorem, we have%
\begin{equation*}
\mu _{F^{M}}(0)=\alpha \beta =M+1.
\end{equation*}

On the other hand, since $M>1$ and both $\lambda _{1}$ and $\lambda _{2}$
are primitive $M$ th roots of unity$,$ the germ $F$ here still satisfies (%
\ref{pf1}). Therefore, by (\ref{pf1}) and by Lemmas \ref{lem2-3} (ii) and %
\ref{lem2-5}, we have%
\begin{equation*}
P_{M}(F,0)=\mu _{F^{M}}(0)-P_{1}(F,0)=M,
\end{equation*}%
and then $\mathcal{O}_{M}(F,0)=1$.
\end{proof}

\begin{description}
\item[\textbf{Part 3}] \textbf{(b3)' implies the existence of }$F$\textbf{.}
\end{description}

\begin{proof}
\textrm{In case (b3)',} we consider the germ $F\in \mathcal{O}(\mathbb{C}%
^{2},0,0)$ given by%
\begin{equation}
F(x_{1},x_{2})=(\lambda _{1}x_{1}+x_{1}^{m_{1}+1}+x_{2}^{d},\lambda
_{2}x_{2}+x_{1}^{m_{1}}x_{2}),  \label{p2+2}
\end{equation}%
where
\begin{equation*}
d=m_{2}/m_{1}=M/m_{1}>1
\end{equation*}%
is a positive integer.

We first assume $m_{1}=1.$ Then $d=m_{2}=M$ and $\lambda _{1}=\lambda
_{2}^{d}=1,$ and then it is easy to show that the $M$ th iteration $%
F^{M}=(F_{1}^{(M)},F_{2}^{(M)})$ has the expression%
\begin{equation*}
\left\{
\begin{array}{l}
F_{1}^{(M)}(x_{1},x_{2})=x_{1}+Mx_{1}^{2}+Mx_{2}^{d}+x_{1}^{2}o(1)+x_{2}^{d}o(1),
\\
F_{2}^{(M)}(x_{1},x_{2})=x_{2}+M\lambda
_{2}^{M-1}x_{1}x_{2}+x_{1}x_{2}o(1)+x_{2}^{d+1}O(1).%
\end{array}%
\right.
\end{equation*}%
Then by Corollary \ref{c1}, considering that $d=M,$ we have that the zero
order $\pi _{F^{M}-id}(0)$ of $F^{M}-id$ at the origin equals $d+2=M+2.$
Thus,%
\begin{equation*}
\mu _{F^{M}}(0)=\pi _{F^{M}-id}(0)=M+2.
\end{equation*}

On the other hand, since $\lambda _{1}=1$ and $\lambda _{2}\neq 1,$ we can
write%
\begin{equation*}
F(x_{1},x_{2})=(x_{1}+x_{1}^{2}+x_{2}o(1),\lambda _{2}x_{2}+x_{2}o(1)),
\end{equation*}%
and then by Corollary \ref{c2}, the zero order $\pi _{F-id}(0)$ of $F-id$ at
the origin equals $2,$ and then
\begin{equation*}
\mu _{F}(0)=P_{1}(F,0)=\pi _{F-id}(0)=2.
\end{equation*}%
Hence, by Lemmas \ref{lem2-5} and \ref{lem2-3} (ii), we have%
\begin{equation*}
P_{M}(F,0)=\mu _{F^{M}}(0)-P_{1}(F,0)=M,
\end{equation*}%
and then $\mathcal{O}_{M}(F,0)=1$.

Now, assume $m_{1}>1.$ Then we have $d=m_{2}/m_{1}>1$ and $\lambda
_{2}^{m_{1}}\neq 1,$ and then, the $m_{1}$ th iteration $%
F^{m_{1}}=(F_{1}^{(m_{1})},F_{2}^{(m_{1})})$ of $F$ given by (\ref{p2+2})
has the expression%
\begin{equation}
\left\{
\begin{array}{l}
F_{1}^{(m_{1})}(x_{1},x_{2})=x_{1}+x_{1}^{m_{1}+1}(c_{1}+o(1))+x_{2}^{d}(c_{1}+o(1)),
\\
F_{2}^{(m_{1})}(x_{1},x_{2})=\lambda _{2}^{m_{1}}x_{2}+m_{1}\lambda
_{2}^{m_{1}-1}x_{1}^{m_{1}}x_{2}(1+o(1))+x_{2}^{d+1}o(1),%
\end{array}%
\right.  \label{p3}
\end{equation}%
where $c_{1}=m_{1}\lambda _{1}^{m_{1}-1}\neq 0.$ Thus, by Corollary \ref{c2}%
, the zero order of the germ $F^{m_{1}}-id$ at the origin equals $m_{1}+1,$
say
\begin{equation}
\mu _{F^{m_{1}}}(0)=m_{1}+1.  \label{p4}
\end{equation}

Now, consider the fixed point index $\mu _{F^{M}}(0)$ of $%
F^{M}=(F_{1}^{(M)},F_{2}^{(M)}),$ the $M$ th iteration of $F$. By \textrm{%
(b3)'} and\ (\ref{p1+2}), $\lambda _{2}^{m_{1}}$ is a primitive $d$ th root
of unity, and then by (\ref{p3}), it is easy to see that the germ (note that
$M=dm_{1}=m_{2})$
\begin{equation*}
F^{M}(x_{1},x_{2})=(F_{1}^{(M)},F_{2}^{(M)})=(F^{m_{1}})^{d}
\end{equation*}%
has the expression%
\begin{equation*}
\left\{
\begin{array}{l}
F_{1}^{(M)}(x_{1},x_{2})=x_{1}+dx_{1}^{m_{1}+1}(c_{1}+o(1))+dx_{2}^{d}(c_{1}+o(1)),
\\
F_{2}^{(M)}(x_{1},x_{2})=x_{2}+M\lambda
_{2}^{M-1}x_{1}^{m_{1}}x_{2}(1+o(1))+x_{2}^{d+1}o(1),%
\end{array}%
\right.
\end{equation*}%
and then the germ $F^{M}-id=(F_{1}^{(M)},F_{2}^{(M)})-id$ has the expression%
\begin{equation*}
\left\{
\begin{array}{l}
F_{1}^{(M)}(x_{1},x_{2})-x_{1}=x_{1}^{m_{1}+1}(dc_{1}+o(1))+x_{2}^{d}(dc_{1}+o(1)),
\\
F_{2}^{(M)}(x_{1},x_{2})-x_{2}=M\lambda
_{2}^{M-1}x_{1}^{m_{1}}x_{2}(1+o(1))+x_{2}^{d+1}o(1),%
\end{array}%
\right.
\end{equation*}%
and then, by Corollary \ref{c1}, we have (note that $dc_{1}\neq 0$ and $%
M\lambda _{2}^{M-1}\neq 0)$
\begin{equation}
\mu _{F^{M}}(0)=\pi _{F^{M}-id}(0)=1+m_{1}+dm_{1}.  \label{p5}
\end{equation}

By Lemmas \ref{lem2-5} and \ref{lem2-3} (ii), We have%
\begin{eqnarray*}
\mu _{F^{M}}(0) &=&P_{M}(F,0)+P_{m_{1}}(F,0)+P_{1}(F,0), \\
\mu _{F^{m_{1}}}(0) &=&P_{m_{1}}(F,0)+P_{1}(F,0).
\end{eqnarray*}%
Thus, by (\ref{p4}) and (\ref{p5}) we have
\begin{equation*}
P_{M}(F,0)=\mu _{F^{M}}(0)-\mu _{F^{m_{1}}}(0)=dm_{1}=M,
\end{equation*}%
and then $\mathcal{O}_{M}(F,0)=1$.
\end{proof}

\begin{description}
\item[Part 4] \textbf{(b4)' implies the existence of }$F$\textbf{.}
\end{description}

\begin{proof}
\textrm{In case (b4)',} the existence of $F$ follows from Proposition \ref%
{c8}. This completes the proof.

Now, we have prove that any case from (b1)' to (b4)' can not occurs. Thus $A$
must satisfy (B) in the main theorem.
\end{proof}

\end{document}